\documentstyle[a4]{article}

\begin{document}

\title{Some local wellposedness results for nonlinear Schr\"odinger equations below $L^2$}

\author{Axel Gr\"unrock \\Fachbereich Mathematik\\ Bergische Universit\"at - Gesamthochschule Wuppertal \\ Gau{\ss}stra{\ss}e 20 \\ D-42097 Wuppertal \\ Germany \\ e-mail Axel.Gruenrock@math.uni-wuppertal.de}

\date{}

\maketitle

\newcommand{\r}{\mbox{${\longrightarrow}$}}
\newcommand{\imp}{\mbox{${\Rightarrow}$}}
\newcommand{\mip}{\mbox{${\Leftarrow}$}}
\newcommand{\iso}{\mbox{${\stackrel{\sim}{\longrightarrow}}$}}
\newcommand{\F}{\mbox{${\cal F}$}}
\newcommand{\Fx}{\mbox{${\cal F}_x$}}
\newcommand{\Ft}{\mbox{${\cal F}_t$}}
\newcommand{\Sn}{\mbox{${\cal S} ({\bf R}^{n+1})$}}
\newcommand{\So}{\mbox{${\cal S} ({\bf R})$}}
\newcommand{\Y}{\mbox{${\cal Y} ({\bf R} \times {\bf T}^n)$}}
\newcommand{\U}{\mbox{$ U_{\phi}$}}
\newcommand{\x}{\mbox{$ X^+_{s,b} $}}
\newcommand{\y}{\mbox{$ Y_{s}$}}
\newcommand{\xm}{\mbox{$ X^-_{s,b} $}}
\newcommand{\ym}{\mbox{$ Y^-_{s}$}}
\newcommand{\xpm}{\mbox{$ X^{\pm}_{s,b} $}}
\newcommand{\ypm}{\mbox{$ Y^{\pm}_{s}$}}
\newcommand{\yy}[1]{\mbox{$ Y_{#1} $}}
\newcommand{\X}[1]{\mbox{$ X_{s,#1}$}}
\newcommand{\Xpm}[1]{\mbox{$ {X^{\pm}}_{s,#1}$}}
\newcommand{\XX}[2]{\mbox{$ X^+_{#1,#2} $}}
\newcommand{\XXm}[2]{\mbox{$ X^-_{#1,#2} $}}
\newcommand{\XXpm}[2]{\mbox{$ X^{\pm}_{#1,#2} $}}
\newcommand{\XXX}[3]{\mbox{$ X_{#1,#2}(#3) $}}
\newcommand{\nfx}{\mbox{$ \| f \| _{ X_{s,b}} $}}
\newcommand{\qfx}{\mbox{$ {\| f \|}^2_{ X_{s,b}} $}}
\newcommand{\nf}[1]{\mbox{$ \| f \| _{ #1} $}}
\newcommand{\qf}[1]{\mbox{$ {\| f \|}^2_{ #1} $}}
\newcommand{\nx}[1]{\mbox{$ \| #1 \| _{ X_{s,b}} $}}
\newcommand{\qx}[1]{\mbox{$ {\| #1 \|}^2_{X_{s,b} } $}}
\newcommand{\n}[2]{\mbox{$ \| #1 \| _{ #2} $}}
\newcommand{\nk}[3]{\mbox{$ \| #1 \| _{ #2} ^{#3}$}}
\newcommand{\q}[2]{\mbox{$ {\| #1 \|}^2_{#2} $}}
\newcommand{\conv}[2]{\mbox{${\mbox{\Huge{$\ast$}}}_{_{\!\!\!\!\!\!\!\!\!{#1}}}^{^{\!\!\!\!\!\!{#2}}}$}}
\pagestyle{plain}
\rule{\textwidth}{0.5pt}

\newtheorem{lemma}{Lemma}[section]
\newtheorem{bem}{Remark}[section]
\newtheorem{bsp}{Example}[section]
\newtheorem{definition}{Definition}[section]
\newtheorem{kor}{Corollary}[section]
\newtheorem{satz}{Theorem}[section]
\newtheorem{prop}{Proposition}[section]

\begin{abstract}
In this paper we prove some local (in time) wellposedness results for nonlinear Schr\"odinger equations
\[u_t - i \Delta u = N(u,\overline{u}) , \hspace{2cm}u(0)=u_0\]
with rough data, that is, the initial value $u_0$ belongs to some Sobolev space of negative index. We obtain positive results for the following nonlinearities and data:
\begin{itemize}
\item $N(u,\overline{u})=\overline{u}^{2}$, \hspace{2,5cm} $u_0 \in H^{s}_x({\bf T}^2)$, \hspace{2.02cm} $s>-\frac{1}{2}$,
\item $N(u,\overline{u})=\overline{u}^{3}$, \hspace{2,6cm}$u_0 \in H^{s}_x({\bf T})$, \hspace{2,27cm}$s>-\frac{1}{3}$,
\item $N(u,\overline{u})=\overline{u}^{2}$, \hspace{2,5cm} $u_0 \in H^{s}_x({\bf T}^3)$, \hspace{1.88cm} $s>-\frac{3}{10}$,
\item $N(u,\overline{u})=u^3$ or $N(u,\overline{u})=\overline{u}^{3}$, 
\hspace{0,25cm}$u_0 \in H^{s}_x({\bf R})$,\hspace{2,13cm} $s>-\frac{5}{12}$,
\item $N(u,\overline{u})=u\overline{u}^{2}$, \hspace{2,3cm} $u_0 \in H^{s}_x({\bf R})$,\hspace{2,26cm} $s>-\frac{2}{5}$,
\item $N(u,\overline{u})=\overline{u}^{4}$, \hspace{2,6cm}$u_0 \in H^{s}_x({\bf T})$ or $u_0 \in H^{s}_x({\bf R})$,\hspace{0,05cm} $s>-\frac{1}{6}$,
\item $N(u,\overline{u})=|u|^{4}$, \hspace{2,41cm}$u_0 \in H^{s}_x({\bf R})$,\hspace{2,25cm} $s>-\frac{1}{8}$,
\item $N(u,\overline{u})=u^4,\,\,\,u^3\overline{u}\,\,\,$ or $u\overline{u}^3$,\,\,\,\hspace{0,54cm}$u_0 \in H^{s}_x({\bf R})$,\hspace{2,25cm} $s>-\frac{1}{6}$.
\end{itemize}
The proof uses the Fourier restriction norm method.
\end{abstract}

\section{Introduction and main results}
The first local (in time) wellposedness results below $L^2$ for the initial value problem for nonlinear Schr\"odinger equations (NLS)
\[u_t - i \Delta u = N(u,\overline{u}) , \hspace{2cm}u(0)=u_0\]
were published in 1996 by Kenig, Ponce and Vega in \cite{KPV96}.
(Here the initial value $u_0$ is assumed to belong to some Sobolev space $H^{s}_x = H^{s}_x({\bf T}^n)$ or $H^{s}_x = H^{s}_x({\bf R}^n)$ with $s<0$.) These authors considered the nonlinearities
\[N_1(u,\overline{u})=u^2,\hspace{1cm}N_2(u,\overline{u})=u\overline{u},\hspace{1cm}N_3(u,\overline{u})=\overline{u}^2\]
in one space dimension. They obtained wellposedness for $N_1$ and $N_3$ under the assumptions
$u_0 \in H^{s}_x({\bf R}),\hspace{0.1cm}s>-\frac{3}{4}$ or $u_0 \in H^{s}_x({\bf T}),\hspace{0.1cm}s>-\frac{1}{2}$
and for $N_2$, provided that
$u_0 \in H^{s}_x({\bf R}),\hspace{0.1cm}s>-\frac{1}{4}$.
This was followed in 1997 by Staffilani's paper \cite{St97}, where wellposedness for NLS with $N=N_3$ and 
$u_0 \in H^{s}_x({\bf R}^2),\hspace{0.1cm}s>-\frac{1}{2}$ 
was shown.

\vspace{0.3cm}

A standard scaling argument suggests that there are even more possible candidates for the nonlinearity to allow local wellposedness below $L^2$: The critical Sobolevexponent for NLS with $N(u,\overline{u})=|u|^{\alpha}u$ obtained by scaling is $s_c=\frac{n}{2}-\frac{2}{\alpha}$. So, for $N_i$, $1\le i \le 3$, there might be local wellposedness for some $s<0$ even for space dimension $n=3$, and in one space dimension also for cubic and quartic nonlinearities positive results seem to be possible.

\vspace{0.3cm}

Recently new theorems concerning this question were presented: In \cite{CDKS} Colliander, Delort, Kenig and Staffilani could prove that in the nonperiodic setting all the results on $N_i$, $1\le i \le 3$, carry over from the one- to the twodimensional case (with the same restrictions on $s$). Concerning the threedimensional nonperiodic case, Tao has shown wellposedness for NLS with the nonlinearities $N_1$ and $N_3$ for $s>-\frac{1}{2}$ and with $N_2$ for $s>-\frac{1}{4}$ (see \cite{T00}, section 11). So concerning the quadratic nonlinearities in the nonperiodic setting the question is meanwhile completely answered.

\vspace{0.3cm}

In this paper the remaining cases are considered, we obtain positive results for the following nonlinearities and data:
\begin{itemize}
\item $N(u,\overline{u})=\overline{u}^{2}$, \hspace{2,5cm} $u_0 \in H^{s}_x({\bf T}^2)$, \hspace{2.02cm} $s>-\frac{1}{2}$,
\item $N(u,\overline{u})=\overline{u}^{3}$, \hspace{2,6cm}$u_0 \in H^{s}_x({\bf T})$, \hspace{2,32cm}$s>-\frac{1}{3}$,
\item $N(u,\overline{u})=\overline{u}^{2}$, \hspace{2,5cm} $u_0 \in H^{s}_x({\bf T}^3)$, \hspace{1.9cm} $s>-\frac{3}{10}$,
\item $N(u,\overline{u})=u^3$ or $N(u,\overline{u})=\overline{u}^{3}$, 
\hspace{0,1cm}$u_0 \in H^{s}_x({\bf R})$,\hspace{2,15cm} $s>-\frac{5}{12}$,
\item $N(u,\overline{u})=u\overline{u}^{2}$, \hspace{2,3cm} $u_0 \in H^{s}_x({\bf R})$,\hspace{2,27cm} $s>-\frac{2}{5}$,
\item $N(u,\overline{u})=\overline{u}^{4}$, \hspace{2,6cm}$u_0 \in H^{s}_x({\bf T})$ or $u_0 \in H^{s}_x({\bf R})$, $s>\!-\frac{1}{6}$,
\item $N(u,\overline{u})=|u|^{4}$, \hspace{2,41cm}$u_0 \in H^{s}_x({\bf R})$,\hspace{2,3cm} $s>-\frac{1}{8}$,
\item $N(u,\overline{u})=u^4,\,\,\,u^3\overline{u}\,\,\,$ or $u\overline{u}^3$,\,\,\,\hspace{0,4cm}$u_0 \in H^{s}_x({\bf R})$,\hspace{2,3cm} $s>-\frac{1}{6}$.
\end{itemize}
To obtain our results, we use the Fourier restriction norm method as it was introduced in \cite{B93} and further developed in \cite{KPV96} and \cite{GTV97}. (In order to concentrate on the crucial multilinear estimates we shall assume this method to be known, for an instructive description thereof we refer to \cite{G96}.) In particular, we will use the function spaces $\xpm = \exp{(\pm it\Delta)} H^b_t(H^s_x)$ equipped with the norms
\begin{eqnarray*}
\n{f}{\xpm}= \n{\exp{(\mp it\Delta)}f}{H^b_t(H^s_x)}=\n{<\xi>^s <\tau \pm |\xi|^2>^b \F f}{L^2_{\xi \tau}} \\
= (\int \mu (d \xi) d \tau <\xi>^{2s}<\tau \pm |\xi|^2>^{2b}|\F f (\xi,\tau)|^2)^{\frac{1}{2}} \hspace{1,5cm}.
\end{eqnarray*}
Here $\F$ denotes the Fourier transform in space and time, $\mu$ is the Lebesgue measure on ${\bf R}^n$ in the nonperiodic respectively the counting measure on ${\bf Z}^n$ in the periodic case, and we use the notation $<x>=(1+|x|^2)^{\frac{1}{2}}$. Observe that $\n{\overline{f}}{\x}=\n{f}{\xm}$. Our proofs rely heavily on the following interpolation property of the $\xpm$-spaces: We have
\[(\XXpm{s_0}{b_0},\XXpm{s_1}{b_1})_{[\theta]} = \xpm \,\,\,,\]
whenever for $\theta \in [0,1]$ it holds that $s=(1-\theta )s_0 + \theta s_1$, $b=(1-\theta )b_0 + \theta b_1$. Here $[\theta]$ denotes the complex interpolation method. Moreover we will make extensive use of the fact that with respect to the inner product on $L^2_{xt}$ the dual space of $\xpm$ is given by $\XXpm{-s}{-b}$. To give a precise formulation of our results, we also need the restriction norm spaces $\xpm (I) = \exp{(\pm it\Delta)} H^b_t(I, H^s_x)$ with norms
\[\n{f}{\xpm (I)}= \inf \{\n{\tilde{f}}{\xpm} : \tilde{f} \in \xpm \,\,\, ,\,\,\,\tilde{f}|_I=f\} .\]
Now our results read as follows:
\begin{satz}Assume
\begin{itemize}
\item[i)] $n=1$, \hspace{3cm} $m=3$, \hspace{3cm}$s>-\frac{1}{3}$ ,\,\,or
\item[ii)] $n=1$, \hspace{3cm} $m=4$,\hspace{3cm} $s>-\frac{1}{6}$ ,\,\,or
\item[iii)] $n=2$, \hspace{3cm} $m=2$,\hspace{3cm} $s>-\frac{1}{2}$ ,\,\,or
\item[iv)] $n=3$,\hspace{3cm}  $m=2$, \hspace{3cm}$s>-\frac{3}{10}$ .
\end{itemize}
Then there exist $b>\frac{1}{2}$ and  $T=T(\n{u_0}{H^s_x({\bf T}^n)})>0$, so that there is a unique solution $u \in \x ([-T,T])$ of the periodic boundary value problem
\[u_t - i \Delta u =  \overline{u}^{m} , \hspace{2cm}u(0)=u_0 \in H^{s}_x({\bf T}^n) .\]
This solution satisfies $u \in C^0([-T,T],H^s_x({\bf T}^n))$ and for any $T' < T$ the mapping \\ $f: H^s_x({\bf T}^n) \r \x ([-T',T'])\,,\,u_0 \mapsto u$ (Data upon solution) is locally Lipschitz continuous.
\end{satz}
The nonlinear estimates leading to this result are contained in Theorems \ref{t1}, \ref{t11} and \ref{t2}, see sections 4 and 5 below. For i) and iii) our results are optimal in the framework of the method and up to the endpoint, in fact there are counterexamples showing that the corresponding multilinear estimates fail for lower values of $s$, see the discussion in section 4. For ii) the scaling argument suggests the optimality of our result. The restriction on $s$ in iv) can possibly be lowered down to $-\frac{1}{2}$, cf. the remark below Thm. \ref{t11}. All the following results are restricted to the onedimensional nonperiodic case:
\begin{satz}Assume 
\begin{itemize}
\item[i)] $s>-\frac{5}{12}$ \hspace{1cm}and \hspace{1cm}$N(u,\overline{u})=u^3$ or $N(u,\overline{u})=\overline{u}^3$, or
\item[ii)] $s>-\frac{2}{5}$ \hspace{2cm}and \hspace{2cm}$N(u,\overline{u})=u\overline{u}^2$.
\end{itemize}
Then there exist $b>\frac{1}{2}$ and $T=T(\n{u_0}{H^s_x})>0$, so that there is a unique solution $u \in \x ([-T,T])$ of the initial value problem
\[u_t - i \partial_x^2 u = N(u,\overline{u})  , \hspace{2cm}u(0)=u_0 \in H^{s}_x ({\bf R}).\]
This solution satisfies $u \in C^0([-T,T],H^s_x)$ and for any $T' < T$ the mapping \\ $f: H^s_x \r \x ([-T',T'])\,,\,u_0 \mapsto u$ (Data upon solution) is locally Lipschitz continuous.
\end{satz}

For the corresponding trilinear estimates see Theorems \ref{t111} and \ref{t1111} (and the remark below) in section 4. We must leave open the question, whether or not the bound on $s$ in the above Theorem can be lowered down to $-\frac{1}{2}$, which is the scaling exponent in this case. This question is closely related to the problem concerning certain trilinear refinements of Strichartz' estimate posed in section 3.

\begin{satz} Let $s >  - \frac{1}{6}$ and $N(u,\overline{u})\in \{u^4, u^3\overline{u},u\overline{u}^3,\overline{u}^4\}$. Then there exist $b>\frac{1}{2}$ and  $T=T(\n{u_0}{H^s_x({\bf R})})>0$, so that there is a unique solution $u \in \x ([-T,T])$ of the initial value problem
\[u_t - i \partial_x^2 u = N(u,\overline{u}), \hspace{2cm}u(0)=u_0 \in H^{s}_x ({\bf R}).\]
This solution satisfies $u \in C^0([-T,T],H^s_x({\bf R}))$ and for any $T' < T$ the mapping $f: H^s_x({\bf R}) \r \x ([-T',T'])\,,\,u_0 \mapsto u$ (Data upon solution) is locally Lipschitz continuous. The same statement holds true for 
$s >  - \frac{1}{8}$ and $N(u,\overline{u})=|u|^4$.
\end{satz}

See Theorems \ref{t2} and \ref{t22} as well as proposition \ref{p51} in section 5
for the crucial nonlinear estimates. The $-\frac{1}{6}$-results should be optimal by scaling, while for the $|u|^4$-nonlinearity the corresponding estimate fails for $s<-\frac{1}{8}$, cf. example \ref{ex53}. Further counterexamples concerned with the periodic case are also given in section 5.
\vspace{0,5cm}

{\bf{Acknowledgement:}} I want to thank Professor Hartmut Pecher for numerous helpful conversations.

\section{Preparatory lemmas}
\subsection{The periodic case}
To prove our results concerning the space-periodic problems, we need the following Strichartz type estimates due to Bourgain:

\begin{lemma}\label{l21} Let $n=1$. Then for all $\epsilon>0$ and $b>\frac{1}{2}$ there exists a constant $c = c(\epsilon,b)$, so that the following estimate holds:
\[\nf{L^6_t({\bf R},L^6_x ({\bf T}))} \leq c \nf{\XX{\epsilon}{b}}\,\,.\]
\end{lemma}

This is essentially Prop. 2.36 in \cite{B93}. For a proof in the form given here, see  \cite{G00}, Lemma 2.2.

\begin{kor}\label{k21} Let $n=1$. Then for all $\epsilon>0$ and $b>\frac{1}{2}$ there exists a constant $c = c(\epsilon,b)$, so that the following estimate holds:
\[\nf{L^8_t({\bf R},L^4_x ({\bf T}))} \leq c \nf{\XX{\epsilon}{b}}\,\,.\]
\end{kor}

Proof: This follows by interpolation between the above lemma and the Sobolev embedding theorem in the time variable.

\begin{lemma}\label{l22} i) Let $n=2$. Then for all $\epsilon>0$ and $b>\frac{1}{2}$ there exists a constant $c = c(\epsilon,b)$, so that the following estimate holds:
\[\nf{L^4_t({\bf R},L^4_x ({\bf T}^2))} \leq c \nf{\XX{\epsilon}{b}}\,\,.\]
ii) Let $n=3$. Then for all $s>\frac{1}{4}$ and $b>\frac{1}{2}$ there exists a constant $c = c(s,b)$, so that the following estimate holds:
\[\nf{L^4_t({\bf R},L^4_x ({\bf T}^3))} \leq c \nf{\XX{s}{b}}\,\,.\]
\end{lemma}

This is essentially the two- respectively the threedimensional case of Prop. 3.6 in \cite{B93}, see also \cite{G00}, Lemma 2.3.

\begin{kor}\label{k22} Let $n=3$. Then for all $s>\frac{1}{5}$ and $b>\frac{9}{20}$ there exists a constant $c = c(s,b)$, so that the following estimate holds:
\[\nf{L^4_t({\bf R},L^{\frac{10}{3}}_x ({\bf T}^3))} \leq c \nf{\XX{s}{b}}\,\,.\]
\end{kor}

Proof: This follows by interpolation between part ii) of the above lemma and the  embedding $\XX{0}{\frac{1}{4}} \subset L^4_t({\bf R},L^2_x ({\bf T}^3))$.

\vspace{0,5cm}

$Remark:$ Because of $\n{f}{L^p_t(L^q_x)}=\n{\overline{f}}{L^p_t(L^q_x)}$ and $\n{f}{\xm}=\n{\overline{f}}{\x}$ the estimates stated in this subsection hold for $\xm$ instead of $\x$. Moreover they are also valid (with $\epsilon =0$) for the corresponding spaces of nonperiodic functions: This is a direct consequence of the Strichartz estimates and \cite{GTV97}, Lemma 2.3.

\subsection{The onedimensional nonperiodic case}

\begin{lemma}\label{l23} Let $n=1$. Then for all $b_0>\frac{1}{2} \ge s \ge 0$,  the following estimates are valid:
\begin{itemize}
\item[i)] $\n{u \overline{v}}{L^2_t(H_x^s)}\le c\n{v}{\XX{0}{b_0}} \n{u}{\XX{0}{b}}$, provided $b>\frac{1}{4}+\frac{s}{2}$,
\item[ii)] $\n{u \overline{v}}{L^p_t(H_x^s)}\le c\n{v}{\XX{0}{b_0}} \n{u}{\XX{0}{b_0}}$, provided $\frac{1}{p}=\frac{1}{4}+\frac{s}{2}$,
\item[iii)] $\n{vw}{\XX{\sigma}{b'}} \le c \n{v}{\XX{\sigma}{b_0}}\n{w}{L^2_t(H_x^{-s-\sigma})}$, provided $\sigma \le 0$, $b'<-\frac{1}{4}-\frac{s}{2}$.
\end{itemize}
\end{lemma}

Proof: We start from the following estimate due to Bekiranov, Ogawa and Ponce
\[\n{u \overline{v}}{L^2_t(\dot{H}_x^{\frac{1}{2}})}\le c\n{u}{\XX{0}{b}} \n{v}{\XX{0}{b}},\hspace{1cm}b>\frac{1}{2}\](see \cite{BOP98}, Lemma 3.2). Combined with\[\n{u \overline{v}}{L^2_{xt}}\le c\n{u}{\XX{0}{b}}\n{v}{\XX{0}{b}},\hspace{1cm}b>\frac{3}{8},\] which follows from Strichartz' estimate, this gives
\begin{equation}\label{21}
\n{u \overline{v}}{L^2_t(H_x^{\frac{1}{2}})}\le c\n{v}{\XX{0}{b_0}} \n{u}{\XX{0}{b}},\hspace{1cm}b_0,b>\frac{1}{2}.
\end{equation}

On the other hand, by H\"older and again by Strichartz' estimate we have
\begin{equation}\label{22}
\n{u \overline{v}}{L^2_{xt}}\le c\n{v}{L^6_{xt}}\n{u}{L^3_{xt}}\le c\n{v}{\XX{0}{b_0}} \n{u}{\XX{0}{b}},\hspace{1cm}b>\frac{1}{4},b_0>\frac{1}{2}.
\end{equation}
Now, by interpolation between (\ref{21}) and (\ref{22}), we obtain part i). To see part ii), we interpolate (\ref{21}) with
\[\n{u \overline{v}}{L_t^4(L_x^2)}\le \n{v}{L_t^8(L_x^4)}\n{u}{L_t^8(L_x^4)} \le c\n{v}{\XX{0}{b_0}} \n{u}{\XX{0}{b_0}}, \hspace{1cm}b_0>\frac{1}{2},\]
which follows from the $L_t^8(L_x^4)$-Strichartz-estimate. Next we dualize part i) to obtain part iii) for $\sigma =0$. For $\sigma < 0$, because of $<\xi_1> \le c <\xi><\xi_2>$, we then have
\[\n{vw}{\XX{\sigma}{b'}} \le c\n{(J^{\sigma}v)(J^{-\sigma}w)}{\XX{0}{b'}} \le c\n{v}{\XX{\sigma}{b_0}}\n{w}{L^2_t(H_x^{-s-\sigma})}.\footnotemark\]
\footnotetext{Here and in the sequel $J^{\sigma}$ ($I^{\sigma}$) denotes the Bessel (Riesz) potential of order $-\sigma$.}
$\hfill \Box$

In order to formulate and prove an analogue for Lemma \ref{l23} in the case of two unbared factors, we introduce some bilinear pseudodifferential operators:

\begin{definition} We define $I_-^s (f,g)$ by its Fourier-transform (in the space variable)
\[\F_x I_-^s (f,g) (\xi) := \int_{\xi_1+\xi_2=\xi}d\xi_1|\xi_1-\xi_2|^s \F_xf(\xi_1)\F_xg(\xi_2).\]
If the expression $|\xi_1-\xi_2|^s$ in the integral is replaced by $<\xi_1-\xi_2>^s$, the corresponding operator will be called $J_-^s (f,g)$. Similarly we define $I_+^s (f,g)$ and $J_+^s (f,g)$ by
\[\F_x I_+^s (f,g) (\xi) := \int_{\xi_1+\xi_2=\xi}d\xi_1|\xi_1+2\xi_2|^s \F_xf(\xi_1)\F_xg(\xi_2).\]
\end{definition}
$Remark$  $(simple \,\,\,properties):$
\begin{itemize}
\item[i)] For functions $u$, $v$ depending on space- and time-variables we have
\[\F I_-^s (u,v) (\xi,\tau) := \int_{\stackrel{\xi_1+\xi_2=\xi}{\tau_1+\tau_2=\tau}}d\xi_1 d\tau_1|\xi_1-\xi_2|^s \F u(\xi_1,\tau_1)\F v(\xi_2,\tau_2)\]
and similar Integrals for the other operators.
\item[ii)] $I_-^s (f,g)$ always coincides with $I_-^s (g,f)$ (and $J_-^s (f,g)$ with $J_-^s (g,f)$), since we can exchange $\xi_1$ and $\xi_2$ in the corresponding integral, while in general we will have $I_+^s (f,g) \neq I_+^s (g,f)$ (and $J_+^s (f,g) \neq J_+^s (g,f)$).
\item[iii)]Fixing $u$ and $s$ we define the linear operators $M$ and $N$ by
\[Mv:=J_-^s (u,v)\hspace{1cm} \mbox{and}\hspace{1cm}Nw:=J_+^s (w,\overline{u}).\]
Then it is easily checked that $M$ and $N$ are formally adjoint with respect to the inner product on $L^2_{xt}$.
\end{itemize}

Now we have the following bilinear Strichartz-type estimate:

\begin{lemma}\label{l24} 
\[\n{I_-^{\frac{1}{2}}(e^{it\partial ^2}u_1,e^{it\partial ^2}u_2)}{L^2_{xt}} \le c \n{u_1}{L^2_{x}}\n{u_2}{L^2_{x}}\]
\end{lemma}

Proof: We will write for short $\hat{u}$ instead of $\F_x u$ and $\int_* d\xi_1$ for $\int_{\xi_1+\xi_2=\xi}d\xi_1$. Then, using Fourier-Plancherel in the space variable we obtain:

\begin{eqnarray*}
&&  \q{I_-^{\frac{1}{2}}(e^{it\partial ^2}u_1,e^{it\partial ^2}u_2)}{L^2_{xt}} \\
&=& c\int d\xi dt |\int_* d\xi_1 |\xi_1-\xi_2|^{\frac{1}{2}} e^{-it(\xi_1^2 + \xi_2^2)}\hat{u}_1(\xi_1)\hat{u}_2(\xi_2)|^2 \\
&=& c\int d\xi dt\int_* d\xi_1 d \eta_1 e^{-it(\xi_1^2 + \xi_2^2- \eta_1^2 - \eta_2^2)}(|\xi_1-\xi_2||\eta_1-\eta_2|)^{\frac{1}{2}}\prod_{i=1}^2 \hat{u_i}(\xi_i) \overline{\hat{u_i}(\eta_i)} \\
&=& c\int d\xi \int_* d\xi_1 d \eta_1 \delta (\eta_1^2 + \eta_2^2 - \xi_1^2 - \xi_2^2)(|\xi_1-\xi_2||\eta_1-\eta_2|)^{\frac{1}{2}}\prod_{i=1}^2 \hat{u_i}(\xi_i) \overline{\hat{u_i}(\eta_i)} \\
&=& c\int d\xi \int_* d\xi_1 d \eta_1 \delta (2(\eta_1^2 - \xi_1^2 + \xi(\xi_1-\eta_1)))(|\xi_1-\xi_2||\eta_1-\eta_2|)^{\frac{1}{2}}\prod_{i=1}^2 \hat{u_i}(\xi_i) \overline{\hat{u_i}(\eta_i)}.
\end{eqnarray*}
Now we use $\delta(g(x)) = \sum_n \frac{1}{|g'(x_n)|} \delta(x-x_n)$, where the sum is taken over all simple zeros of $g$, in our case:
\[g(x)= 2(x^2 + \xi (\xi_1-x)-\xi_1^2)\]
with the zeros $x_1 =\xi_1$ and $x_2 =\xi-\xi_1$, hence $g'(x_1)=2(2\xi_1-\xi)$ respectively $g'(x_2)=2(\xi-2\xi_1)$.
So the last expression is equal to
\begin{eqnarray*}
&&c\int d\xi \int_* d\xi_1 d \eta_1\frac{1}{|2\xi_1-\xi|}\delta(\eta_1-\xi_1)(|\xi_1-\xi_2||\eta_1-\eta_2|)^{\frac{1}{2}}\prod_{i=1}^2 \hat{u_i}(\xi_i) \overline{\hat{u_i}(\eta_i)}\\
&+& c\int d\xi \int_* d\xi_1 d \eta_1\frac{1}{|2\xi_1-\xi|}\delta(\eta_1-(\xi-\xi_1))(|\xi_1-\xi_2||\eta_1-\eta_2|)^{\frac{1}{2}}\prod_{i=1}^2 \hat{u_i}(\xi_i) \overline{\hat{u_i}(\eta_i)}\\
&=&c\int d\xi \int_* d\xi_1 \prod_{i=1}^2|\hat{u_i}(\xi_i)|^2 + c\int d\xi \int_* d\xi_1\hat{u}_1(\xi_1)\overline{\hat{u}_1}(\xi_2)\hat{u}_2(\xi_2)\overline{\hat{u}_2}(\xi_1)\\
&\le& c (\prod_{i=1}^2 \q{u_i}{L^2_x} + \q{\hat{u}_1\hat{u}_2}{L^1_{\xi}}) \le c \prod_{i=1}^2 \q{u_i}{L^2_x}.
\end{eqnarray*}
$\hfill \Box$

\vspace{0.5cm}

\begin{kor}\label{k23} Let $b_0 > \frac{1}{2}$ and $0 \le s \le \frac{1}{2}$. Then the following estimates hold true:
\begin{itemize}
\item[i)] $\n{J^s_- (u,v)}{L^2_{xt}}\le c \n{u}{\XX{0}{b_0}}\n{v}{\XX{0}{b}}$, provided $b>\frac{1}{4}+\frac{s}{2}$,
\item[ii)] $\n{J^s_+ (v,\overline{u})}{\XX{0}{b'}}\le c \n{u}{\XX{0}{b_0}}\n{v}{L^2_{xt}}$, provided $b'>-\frac{1}{4}-\frac{s}{2}$.
\end{itemize}
\end{kor}

$Remark:$ In i) we may replace $J^s_- (u,v)$ by $J^s_- (\overline{u},\overline{v})$, in fact a short computation shows that $J^s_- (\overline{u},\overline{v})=\overline{J^s_- (u,v)}$.

\vspace{0.5cm}

Proof: Arguing as in the proof of Lemma 2.3 in \cite{GTV97}, we obtain from the above Lemma
\[\n{I^{\frac{1}{2}}_-(u,v)}{L^2_{xt}} \le c \n{u}{\XX{0}{b_0}}\n{v}{\XX{0}{b}},\hspace{1cm}b,b_0 >\frac{1}{2}.\]
Combining this with
\[\n{uv}{L^2_{xt}}\le \n{u}{L^6_{xt}}\n{v}{L^3_{xt}}\le c\n{u}{\XX{0}{b_0}} \n{v}{\XX{0}{b}},\hspace{1cm}b>\frac{1}{4},b_0>\frac{1}{2},\]
we obtain i) for $s=\frac{1}{2}$ and $s=0$.

To see i) for $0<s<\frac{1}{2}$, $b>\frac{1}{4}+\frac{s}{2}$, we write $w=\Lambda^b v$, where $\Lambda^b$ is defined by $\F \Lambda^b v(\xi,\tau)=<\tau+\xi^2>^b \F v(\xi,\tau)$. Then we have to show that
\begin{equation}\label{23}
\n{J^s_- (u,\Lambda^{-b}w)}{L^2_{xt}}\le c \n{u}{\XX{0}{b_0}}\n{w}{L^2_{xt}},
\end{equation}
where
\[\n{J^s_- (u,\Lambda^{-b}w)}{L^2_{xt}}=\n{\int_{\stackrel{\tau_1+\tau_2=\tau}{\xi_1+\xi_2=\xi}}<\!\!\xi_1-\xi_2\!\!>^s\F u(\xi_1,\tau_1)<\!\!\tau_2+\xi_2^2\!\!>^{-b}\F w(\xi_2,\tau_2)}{L^2_{\xi \tau}}.\]
Notice that, by the preceding, (\ref{23}) is already known in the limiting cases $(s,b)=(0,\frac{1}{4}+\epsilon))$ and $(s,b)=(\frac{1}{2},\frac{1}{2}+\epsilon)$, $\epsilon>0$. Choosing $\epsilon=b-\frac{1}{4}-\frac{s}{2}$ we have
\[<\!\!\xi_1-\xi_2\!\!>^s <\!\!\tau_2+\xi_2^2\!\!>^{-b} \le <\!\!\tau_2+\xi_2^2\!\!>^{-\frac{1}{4}-\epsilon} +      <\!\!\xi_1-\xi_2\!\!>^{\frac{1}{2}}<\!\!\tau_2+\xi_2^2\!\!>^{-\frac{1}{2}-\epsilon}\]
and hence
\[\n{J^s_- (u,\Lambda^{-b}w)}{L^2_{xt}}\le \n{ u(\Lambda^{-\frac{1}{4}-\epsilon}w)}{L^2_{xt}}+\n{J^{\frac{1}{2}}_- (u,\Lambda^{-\frac{1}{2}-\epsilon}w)}{L^2_{xt}}\le c\n{u}{\XX{0}{b_0}}\n{w}{L^2_{xt}}.\]

Finally, ii) follows from i) by duality (cf. part iii) of the remark on simple properties of $J^s_-$).
$\hfill \Box$

\section{Trilinear refinements of the onedimensional $L^6$-Strichartz-estimate in the nonperiodic case}

In \cite{B98} Bourgain showed the following bilinear refinement of the $L^4_{xt}$-Strichartz-estimate in two space dimensions 
\[\n{u_1u_2}{L^2_t(H^s_x)}\le c \n{u_1}{\XX{s+\epsilon}{b}}\n{u_2}{\XX{0}{b}},\]
provided $0 \le s < \frac{1}{2}<b$, $\epsilon>0$. The exponent in the onedimensional Strichartz estimate is $6$, so the question for trilinear refinements of this estimate comes up naturally. In this section we shall give a partial answer to this question, starting with the following fairly easy application of Kato's smoothing effect:

\begin{lemma}\label{l31} Let $0\le s\le \frac{1}{4}$, $b>\frac{1}{2}$. Then the estimate
\[\n{u_1 u_2 u_3}{L^2_{xt}} \le c \n{u_1}{\XX{s}{b}}\n{u_2}{\XX{-s}{b}}\n{u_3}{\XX{0}{b}}\]
holds true.
\end{lemma}
 
Proof: For $s=0$ this follows from standard Strichartz' estimate, for $s=\frac{1}{4}$ we argue as follows: Interpolation between the $L^6$-estimate and the Kato smoothing effect
\[\n{e^{it\partial ^2}u_0}{L_x^{\infty}(L_t^2)}\le c \n{u_0}{\dot{H}_x^{-\frac{1}{2}}}\]
(see Thm. 4.1 in \cite{KPV91}) with $\theta =\frac{1}{2}$ yields
\[\n{I^{\frac{1}{4}} e^{it\partial ^2}u_0}{L_x^{12}(L_t^3)}\le c\n{u_0}{L^2_{x}}.\]
Now Lemma 2.3 in \cite{GTV97} gives
\begin{equation}\label{31}
\n{I^{\frac{1}{4}}u}{L_x^{12}(L_t^3)} \le c \n{u}{\XX{0}{b}}, \,\,\,\,b>\frac{1}{2}.
\end{equation}
On the other hand by Thm. 2.5 in \cite{KPV91} we get
\[\n{e^{it\partial ^2}u_0}{L_x^4(L_t^{\infty})}\le c \n{u_0}{\dot{H}_x^{\frac{1}{4}}}\]
and thus
\begin{equation}\label{32}
\n{u}{L_x^4(L_t^{\infty})}\le c \n{I^{\frac{1}{4}}u}{\XX{0}{b}}\le c \n{u}{\XX{\frac{1}{4}}{b}}, \,\,\,\,b>\frac{1}{2}.
\end{equation}
Using the projections $p$ and $P$ defined by $p=\F^{-1}\chi_{\{|\xi| \le 1\}}\F$ and $P= Id-p$, we now have
\[\n{u_1u_2}{L^3_{xt}}\le \n{u_1pu_2}{L^3_{xt}}+\n{u_1Pu_2}{L^3_{xt}}=:N_1+N_2\]
with
\[N_1 \le \n{u_1}{L^6_{xt}}\n{pu_2}{L^6_{xt}}\le c\n{u_1}{\XX{0}{b}}\n{pu_2}{\XX{0}{b}}\le c\n{u_1}{\XX{\frac{1}{4}}{b}}\n{u_2}{\XX{-\frac{1}{4}}{b}}.\]
For $N_2$ we use (\ref{31}) and (\ref{32}) to obtain
\begin{eqnarray*}
\n{u_1Pu_2}{L^3_{xt}}&\le &\n{u_1}{L_x^4(L_t^{\infty})}\n{Pu_2}{L_x^{12}(L_t^3)}\\
&\le & c \n{u_1}{\XX{\frac{1}{4}}{b}}\n{I^{-\frac{1}{4}}Pu_2}{\XX{0}{b}}\le c\n{u_1}{\XX{\frac{1}{4}}{b}}\n{u_2}{\XX{-\frac{1}{4}}{b}}.
\end{eqnarray*}
Now, using H\"older and standard Strichartz again, from this we obtain the claim for $s=\frac{1}{4}$. For $0<s<\frac{1}{4}$ the result then follows by multilinear interpolation, see Thm. 4.4.1 in \cite{BL}.
$\hfill \Box$

\vspace{0,3cm}

{\bf Problem:} Does the above estimate hold for $\frac{1}{4}<s<\frac{1}{2}$ ?

\vspace{0,3cm}
\begin{kor}\label{k31} Assume $0\le s\le \frac{1}{4}$ and $b>\frac{1}{2}$. Let $\tilde{u}$ denote $u$ or $\overline{u}$. Then the following estimates are valid:
\begin{itemize}
\item[i)] $\n{\tilde{u}_1 \tilde{u}_2 \tilde{u}_3}{L^2_{xt}} \le c \n{u_1}{\XX{s}{b}}\n{u_2}{\XX{-s}{b}}\n{u_3}{\XX{0}{b}}$,
\item[ii)] $\n{\tilde{u}_1 \tilde{u}_2 \tilde{u}_3}{\XX{-s}{-b}} \le c \n{u_1}{L^2_{xt}}\n{u_2}{\XX{-s}{b}}\n{u_3}{\XX{0}{b}}$,
\item[iii)] $\n{\tilde{u}_1 \tilde{u}_2 \tilde{u}_3}{L^2_t (H^s_x)} \le c \n{u_1}{\XX{s}{b}}\n{u_2}{\XX{0}{b}}\n{u_3}{\XX{0}{b}}$,
\item[iv)] $\n{\tilde{u}_1 \tilde{u}_2 \tilde{u}_3}{\XX{-s}{-b}} \le c \n{u_1}{L^2_t (H^{-s}_x)}\n{u_2}{\XX{0}{b}}\n{u_3}{\XX{0}{b}}$.
\end{itemize}
\end{kor}

Proof: Clearly, in $\n{u_1 u_2 u_3}{L^2_{xt}}$ any factor $u_i$ may be replaced by $\overline{u}_i$. This gives i). From this we obtain ii) by duality. Writing $<\xi> \le <\xi_1>+<\xi_2>+<\xi_3>$ and applying i) twice (plus standard Strichartz), part iii) can be seen. Dualizing again, part iv) follows.
$\hfill \Box$

\vspace{0,3cm}

In some cases, using the bilinear estimates of the previous section, we can prove better $L^2_t(H^s_x)$-estimates:

\begin{lemma}\label{l32}
\begin{itemize}
\item[i)] For $|s|<\frac{1}{2}<b$ the following estimate holds:
\[\n{u_1 \overline{u}_2u_3}{L_t^2(H_x^s)} \le c\n{u_1}{\XX{0}{b}}\n{u_2}{\XX{0}{b}}\n{u_3}{\XX{s}{b}}\]
\item[ii)] For $-\frac{1}{2}<s \le 0$, $b>\frac{1}{2}$ the following is valid:
\[\n{u_1 \overline{u}_2u_3}{L_t^2(H^{s} _x)} \le c\n{u_1}{\XX{0}{b}}\n{u_2}{\XX{s}{b}}\n{u_3}{\XX{0}{b}}\]
\end{itemize}
\end{lemma}

\vspace{0,3cm}

$Remark:$ Using multilinear interpolation (Thm. 4.4.1 in \cite{BL}) we obtain
\[\n{u_1 \overline{u}_2u_3}{L_t^2(H^{s} _x)} \le c\n{u_1}{\XX{s_1}{b}}\n{u_2}{\XX{s_2}{b}}\n{u_3}{\XX{s_3}{b}},\]
provided $-\frac{1}{2}<s \le 0$, $b>\frac{1}{2}$, $s_{1,2,3} \le 0$ and $s_1+s_2+s_3=s$. Moreover, we may replace $u_1 \overline{u}_2u_3$ on the left hand side by $\overline{u}_1 u_2 \overline{u}_3$.

\vspace{0,3cm}

Proof: First we show i) for $s > 0$. From $<\xi>\le c (<\xi_1+\xi_2>+<\xi_3>)$ it follows that
\[\n{u_1 \overline{u}_2u_3}{L_t^2(H_x^s)} \le c\n{J^s(u_1 \overline{u}_2)u_3}{L_{xt}^2} + \n{u_1 \overline{u}_2J^su_3}{L_{xt}^2} =: c (N_1+N_2).\]
Using the standard $L^6_{xt}$-Strichartz-estimate we see that $N_2$ is bounded by the right hand side of i). For $N_1$ we have with $s=\frac{1}{p}$, $\frac{1}{2}-s=\frac{1}{q}$ ($\imp H^s \subset L^q,\,\,\,H^{\frac{1}{2}} \subset H^{s,p}$):
\begin{eqnarray*}
N_1 &\le& c \n{J^s(u_1 \overline{u}_2)}{L_t^2(L_x^p)}\n{u_3}{L_t^{\infty}(L_x^q)}\\
&\le& c \n{u_1 \overline{u}_2}{L_t^2(H^{\frac{1}{2}}_x)}\n{u_3}{L_t^{\infty}(H^s_x)}\\
&\le& c \n{u_1}{\XX{0}{b}}\n{u_2}{\XX{0}{b}}\n{u_3}{\XX{s}{b}}
\end{eqnarray*}
by Lemma \ref{l23}, i), and the Sobolev embedding in the time variable.

Next we consider i) for $s<0$. Writing $<\xi_3>\le c(<\xi>+<\xi_1 +\xi_2>)$, we obtain
\[\n{u_1 \overline{u}_2u_3}{L_t^2(H_x^s)} \le c\n{u_1 \overline{u}_2J^s u_3}{L_{xt}^2} + \n{J^{-s}(u_1 \overline{u}_2)J^su_3}{L_t^2(H^s_x)} =: c (N_1+N_2).\]
To estimate $N_1$ we use again the standard $L^6_{xt}$-Strichartz estimate. For $N_2$ we use the embedding $L^q \subset H^s,\,\,s-\frac{1}{2}=-\frac{1}{q}$ and H\"older's inequality:
\begin{eqnarray*}
N_2 &\le& c \n{J^{-s}(u_1 \overline{u}_2)J^su_3}{L_t^2(L^q_x)}\\
&\le& c \n{J^{-s}(u_1 \overline{u}_2)}{L_t^2(L^p_x)}\n{u_3}{L_t^{\infty}(H^s_x)},
\end{eqnarray*}
where $\frac{1}{q}=\frac{1}{2}+\frac{1}{p}$. The second factor is bounded by $c\n{u_3}{\XX{s}{b}}$ because of Sobolev's embedding Theorem in the time variable. For the first factor we use the embedding $H^{\frac{1}{2}}\subset H^{-s,p}$ (observe that $s=-\frac{1}{p}$) and again Lemma \ref{l23}, i).

We conclude the proof by showing ii): Here we have $\xi =(\xi_1+\xi_2)+(\xi_3+\xi_2)-\xi_2$ respectively $<\xi_2>\le c(<\xi>+<\xi_1 +\xi_2>+<\xi_3 +\xi_2>)$ and thus
\[\n{u_1 \overline{u}_2u_3}{L_t^2(H_x^s)} \le c(N_1 +N_2 +N_3)\]
with
\[N_1 =\n{u_1 (J^s\overline{u}_2) u_3}{L_{xt}^2} \le c\n{u_1}{\XX{0}{b}}\n{u_2}{\XX{s}{b}}\n{u_3}{\XX{0}{b}}\]
(by standard Strichartz) and
\[N_2= \n{J^{-s}(u_1 J^s\overline{u}_2)u_3}{L_t^2(H^s_x)},\,\,\,\,\,N_3=\n{u_1 J^{-s}((J^s\overline{u}_2)u_3)}{L_t^2(H^s_x)}.\]
By symmetry between $u_1$ and $u_3$ it is now sufficient to estimate $N_2$: Using the embedding $L^q \subset H^s,\,\,\,s-\frac{1}{2}=-\frac{1}{q}$, H\"older's inequality and the embedding $H^{\frac{1}{2}}\subset H^{-s,p},\,\,\,-s=\frac{1}{p}$ we obtain
\begin{eqnarray*}
N_2 &\le& c \n{J^{-s}(u_1 J^s\overline{u}_2)u_3}{L_t^2(L^q_x)}\\
&\le & c \n{J^{-s}(u_1 J^s\overline{u}_2)}{L_t^2(L^p_x)}\n{u_3}{L_t^{\infty}(L^2_x)}\\
&\le & c \n{J^{\frac{1}{2}}(u_1 J^s\overline{u}_2)}{L^2_{xt}}\n{u_3}{L_t^{\infty}(L^2_x)}.
\end{eqnarray*}
Again, Lemma \ref{l23}, i), and the Sobolev embedding in $t$ give the desired bound.
$\hfill \Box$

\begin{lemma}\label{l33} 
For $-\frac{1}{2}<s \le 0$, $b>\frac{1}{2}$ the following holds true:
\[\n{u_1 u_2 u_3}{L_t^2(H^{s} _x)} \le c\n{u_1}{\XX{s}{b}}\n{u_2}{\XX{0}{b}}\n{u_3}{\XX{0}{b}}\]
\end{lemma}

\vspace{0,3cm}

$Remark:$ Again we may use multilinear interpolation to get
\[\n{u_1 u_2u_3}{L_t^2(H^{s} _x)} \le c\n{u_1}{\XX{s_1}{b}}\n{u_2}{\XX{s_2}{b}}\n{u_3}{\XX{s_3}{b}}\]
for $-\frac{1}{2}<s \le 0$, $b>\frac{1}{2}$, $s_{1,2,3} \le 0$ and $s_1+s_2+s_3=s$. The same holds true with $u_1 u_2u_3$ replaced by $\overline{u}_1 \overline{u}_2 \overline{u}_3$.
\vspace{0,3cm}

Proof: It is easily checked that for $\rho, \lambda \ge 0$ the inequality
\[<\xi_1>^{\rho} \le c(<\xi>^{\rho}+\frac{<\xi_1-\xi_2>^{\rho+\lambda}}{<\xi_1+\xi_2>^{\lambda}}+\frac{<\xi_1-\xi_3>^{\rho+\lambda}}{<\xi_1+\xi_3>^{\lambda}})\]
is valid, if $\xi=\xi_1+\xi_2+\xi_3$. Choosing $\rho = -s$ and $\lambda = s+\frac{1}{2}$ it follows, that
\[\n{u_1 u_2 u_3}{L_t^2(H^{s} _x)} \le c (N_1+N_2+N_3),\]
where
\[N_1=\n{(J^su_1) u_2 u_3}{L^2_{xt}} \le c \n{u_1}{\XX{s}{b}}\n{u_2}{\XX{0}{b}}\n{u_3}{\XX{0}{b}}\]
(by standard Strichartz) and
\[N_2= \n{(J^{-\lambda}J_-^{\frac{1}{2}}(J^su_1, u_2))u_3}{L_t^2(H^s_x)},\,\,\,\,\,N_3=\n{(J^{-\lambda}J_-^{\frac{1}{2}}(J^su_1, u_3))u_2}{L_t^2(H^s_x)}.\]
Now, by symmetry between $u_2$ and $u_3$, it is sufficient to estimate $N_2$. Using the embedding $L^q \subset H^s$, ($s-\frac{1}{2}=-\frac{1}{q}$) and H\"older we get
\begin{eqnarray*}
N_2 &\le& c \n{J^{-\lambda}J_-^{\frac{1}{2}}(J^su_1, u_2)u_3}{L_t^2(L^q_x)}\\
&\le & c \n{J^{-\lambda}J_-^{\frac{1}{2}}(J^su_1, u_2)}{L_t^2(L^p_x)}\n{u_3}{L_t^{\infty}(L^2_x)}
\end{eqnarray*}
with $\frac{1}{q}=\frac{1}{2}+\frac{1}{p}$. The second factor is bounded by $c\n{u_3}{\XX{0}{b}}$. For the first factor we observe that $L^2 \subset H^{-\lambda,p}$, so it can be estimated by
\[\n{J_-^{\frac{1}{2}}(J^su_1, u_2)}{L^2_{xt}}\le c \n{u_1}{\XX{s}{b}}\n{u_2}{\XX{0}{b}},\]
where in the last step we have used Corollary \ref{k23}, i).
$\hfill \Box$

\section{Estimates on quadratic and cubic nonlinearities}
\begin{satz}\label{t1} Let $n=1, m=3$ or $n=2, m=2$. Assume $ 0 \geq s > - \frac{1}{m}$ and $ - \frac{1}{2} < b' < \frac{ms}{2}$. Then in the periodic and nonperiodic case for all $b>\frac{1}{2}$ the estimate 
\[\n{ \prod_{i=1}^m\overline{u}_i}{\XX{0}{b'}} \leq c \prod_{i=1}^m \n{u_i}{\XX{s}{b}}\]
holds true.
\end{satz}

Proof: Defining $f_i(\xi, \tau)= <\tau - |\xi|^2>^{b}<\xi>^s \F \overline{u}_i(\xi, \tau)$, $1 \leq i \leq m$, we have 
\[\n{ \prod_{i=1}^m \overline{u}_i}{\XX{0}{b'}}=c \n{<\!\!\tau + |\xi|^2\!\!>^{b'} \int d \nu  \prod_{i=1}^m <\!\!\tau_i - |\xi_i|^2\!\!>^{-b} <\xi_i>^{-s}f_i(\xi_i,\tau_i)}{L^2_{\xi,\tau}},\]
where $d \nu = \mu(d\xi_1..d\xi_{m-1}) d\tau_1.. d \tau_{m-1}$ and $\sum_{i=1}^m (\xi_i,\tau_i) = (\xi, \tau)$. \footnote{In the sequel we shall make repeated use of this convention.} Because of
\[\tau + |\xi|^2 - \sum_{i=1}^m\tau_i - |\xi_i|^2 =|\xi|^2 + \sum_{i=1}^m |\xi_i|^2\]
there is the inequality
\begin{eqnarray}\label{41}
<\xi>^2 + \sum_{i=1}^m <\xi_i>^2 & \le & <\tau + |\xi|^2 > + \sum_{i=1}^m <\tau_i - |\xi_i|^2> \nonumber \\
& \le & c (<\tau + |\xi|^2 > + \sum_{i=1}^m <\tau_i - |\xi_i|^2> \chi_{A_i}),
\end{eqnarray}
where in $A_i$ we have $<\tau_i - |\xi_i|^2> \,\,\ge \,\, <\tau + |\xi|^2 >$. Since $b' < \frac{ms}{2}$ is assumed, it follows
\[<\xi>^{\epsilon}\prod_{i=1}^m<\xi_i>^{-s+\epsilon} \leq c (<\tau + |\xi|^2 >^{-b'} + \sum_{i=1}^m <\tau_i - |\xi_i|^2>^{-b'} \chi_{A_i})\]
for some $\epsilon > 0$. From this we conclude that
\[\n{ \prod_{i=1}^m \overline{u}_i}{\XX{0}{b'}} \le c \sum_{j=0}^m \n{I_j}{L^2_{\xi,\tau}},\]
with
\[I_0(\xi,\tau)=<\xi>^{-\epsilon}\int d \nu  \prod_{i=1}^m <\tau_i - |\xi_i|^2>^{-b} <\xi_i>^{-\epsilon}f_i(\xi_i,\tau_i)\]
and, for $1 \le j \le m$,
\begin{eqnarray*}
I_j(\xi,\tau)\!\!\!\! & = & \!\!\!\!<\!\!\xi\!\!>\!\!^{-\epsilon}<\!\!\tau \!\!+ \!\!|\xi|^2 \!\!>^{b'}\!\!\int\!\! d \nu <\!\!\tau_j \!\!- \!\!|\xi_j|^2 \!\! >\!\!\!\!^{-b'} \prod_{i=1}^m <\!\!\tau_i \!\!- \!\!|\xi_i|^2\!\!>\!\!\!\!^{-b} <\!\!\xi_i\!\!>\!\!^{-\epsilon}f_i(\xi_i,\tau_i)\chi_{A_i} \\
\!\!\!\!& \le &\!\!\!\! <\!\!\xi\!\!>^{-\epsilon}<\!\!\tau \!\!+ \!\!|\xi|^2 \!\!>^{-b}\!\!\int \!\! d \nu <\!\!\tau_j\!\! - \!\!|\xi_j|^2\!\! >^{b} \prod_{i=1}^m <\!\!\tau_i\!\! -\!\! |\xi_i|^2\!\!>^{-b} <\!\!\xi_i\!\!>^{-\epsilon} f_i(\xi_i,\tau_i).
\end{eqnarray*}
To estimate $I_0$ we use H\"olders inequality and Lemma \ref{l21} respectively \ref{l22}:
\begin{eqnarray*}
\n{I_0}{L^2_{\xi,\tau}} & \le & \n{\int d \nu  \prod_{i=1}^m <\tau_i - |\xi_i|^2>^{-b} <\xi_i>^{-\epsilon}f_i(\xi_i,\tau_i)}{L^2_{\xi,\tau}}\\
&=& c \n{\prod_{i=1}^m J^{s-\epsilon} \overline{u}_i}{L^2_{x,t}}
 \le c\prod_{i=1}^m\n{ J^{s-\epsilon} \overline{u}_i}{L^{2m}_{x,t}}\\
& \le &c\prod_{i=1}^m\n{ J^{s} \overline{u}_i}{\XXm{0}{b}}=c\prod_{i=1}^m\n{\overline{u}_i}{\XXm{s}{b}}.
\end{eqnarray*}
To estimate $I_j$, $1\le j \le m$, we define  $p=2m$  and  $p'$  by $\frac{1}{p}+\frac{1}{p'}=1$. Then we use the dual versions of Lemma \ref{l21} respectively \ref{l22}, H\"olders inequality and the Lemmas themselves to obtain:
\begin{eqnarray*}
\n{I_j}{L^2_{\xi,\tau}} 
& \le & c\n{(\prod_{\stackrel{i=1}{i \ne j}}^m J^{s-\epsilon} \overline{u}_i)(J^{-\epsilon}\F^{-1}f_j)}{\XX{-\epsilon}{-b}}\\
& \le & c\n{(\prod_{\stackrel{i=1}{i \ne j}}^m J^{s-\epsilon} \overline{u}_i)(J^{-\epsilon}\F^{-1}f_j)}{L^{p'}_{x,t}}\\
& \le & c\n{J^{-\epsilon}\F^{-1}f_j}{L^2_{x,t}}\prod_{\stackrel{i=1}{i \ne j}}^m\n{ J^{s-\epsilon} \overline{u}_i}{L^{p}_{x,t}}\\
& \le & c\n{f_j}{L^2_{\xi,\tau}}\prod_{\stackrel{i=1}{i \ne j}}^m\n{ J^{s} \overline{u}_i}{\XXm{0}{b}} = c \prod_{i=1}^m\n{\overline{u}_i}{\XXm{s}{b}}
\end{eqnarray*}
$\hfill \Box$

\vspace{0,5cm}

$Remark:$ The above theorem with $n=m=2$ can be inserted into the proof of Theorem 2.5 in \cite{St97}, thus showing that the statement of that theorem also holds in the periodic case. So we can answer this question left open in \cite{St97} affirmatively (cf. the remark on top of p. 81 in \cite{St97}). Moreover it is a straightforward application of Sobolev's embedding theorem, to prove the complementary estimate
\[\n{\prod_{i=1}^m u_i}{\XX{0}{-b'}} \le c \prod_{i=1}^m \n{u_i}{\XXpm{1}{b}},\,\,\,-b'< \frac{1}{2}<b,\]
$m \ge 2$ arbitrary. (Observe that Thm. \ref{t1} holds with $s=0$ on the left hand side.) So we obtain the bound
\[\n{u(t)}{H^s_x ({\bf T}^2)} \le c <t>^{s-1+ \epsilon} ,\,\,\,s>1,\,\,\,\epsilon>0,\]
whenever $u$ is a global solution of
\[iu_t + \Delta u + \lambda |u|^{2l} u =0, \hspace{1cm}u(0)=u_0 \in H^s_x ({\bf T}^2),\]
($l \in {\bf N}$) and $\n{u(t)}{H^1_x ({\bf T}^2)}$ is controlled by the conserved energy.

\begin{satz}\label{t11} Let $n=3$ and assume $ 0 \geq s > - \frac{3}{10}$, $ - \frac{1}{2} < b' < \frac{s}{2} - \frac{7}{20}$ and $b>\frac{1}{2}$. Then in the periodic case the estimate 
\[\n{ \prod_{i=1}^2\overline{u}_i}{\XX{s}{b'}} \leq c \prod_{i=1}^2 \n{u_i}{\XX{s}{b}}\]
holds true.
\end{satz}

Proof: Writing $f_i(\xi, \tau)= <\tau - |\xi|^2>^{b}<\xi>^s \F \overline{u}_i(\xi, \tau)$, $1 \leq i \leq 2$, we have
\begin{eqnarray*}
\n{ \prod_{i=1}^2\overline{u}_i}{\XX{s}{b'}}  \hspace{4cm}\\
= c \n{<\xi>^s<\tau + |\xi|^2>^{b'} \int d \nu  \prod_{i=1}^2 <\tau_i - |\xi_i|^2>^{-b} <\xi_i>^{-s}f_i(\xi_i,\tau_i)}{L^2_{\xi,\tau}}.
\end{eqnarray*}
By the expressions $<\tau + |\xi|^2>$ and $<\tau_i - |\xi_i|^2>$, $i=1,2$, the quantity \\ $<\xi>^2 +<\xi_1>^2+<\xi_2>^2$ can be controlled. So we split the domain of integration into $A_0 +A_1 +A_2$, where in $A_0$ we have \\ $<\tau + |\xi|^2> = \max{(<\tau + |\xi|^2>, <\tau_{1,2} - |\xi_{1,2}|^2>)}$ and in $A_j$, $j=1,2$, it should hold that $<\tau_j - |\xi_j|^2> = \max{(<\tau + |\xi|^2>, <\tau_{1,2} - |\xi_{1,2}|^2>)}$. First we consider the region $A_0$: Here we use that for $\epsilon > 0$ sufficiently small
\[<\xi>^{\frac{3}{10} + s} \prod_{i=1}^2<\xi_i>^{-s + \frac{1}{5}+ \epsilon}\le c <\tau + |\xi|^2>^{-b'}.\]
This gives the upper bound
\begin{eqnarray*}
\n{<\xi>^{-\frac{3}{10}}\int d \nu  \prod_{i=1}^2 <\tau_i - |\xi_i|^2>^{-b} <\xi_i>^{-\frac{1}{5}- \epsilon}f_i(\xi_i,\tau_i)}{L^2_{\xi,\tau}} \\
=c \n{\prod_{i=1}^2 J^{s-\frac{1}{5}- \epsilon}\overline{u}_i}{L_t^2(H_x^{-\frac{3}{10}})}.\hspace{3cm}
\end{eqnarray*}
Now, using the embedding $L_x^q \subset H_x^{-\frac{3}{10}}$, $\frac{1}{q}=\frac{3}{5}$, H\"older's inequality and Corollary \ref{k22}, we get the following chain of inequalities:
\begin{eqnarray*}
\n{\prod_{i=1}^2 J^{s-\frac{1}{5}- \epsilon}\overline{u}_i}{L_t^2(H_x^{-\frac{3}{10}})} &\le & c 
\n{\prod_{i=1}^2 J^{s-\frac{1}{5}- \epsilon}\overline{u}_i}{L_t^2(L_x^q)} \\
&\le & c \n{J^{s-\frac{1}{5}- \epsilon}u_1}{L_t^4(L_x^{2q})}\n{J^{s-\frac{1}{5}- \epsilon}u_2}{L_t^4(L_x^{2q})} \\
& \le & c \prod_{i=1}^2 \n{u_i}{\XX{s}{b}}.
\end{eqnarray*}

Now, by symmetry, it only remains to show the estimate for the region $A_1$: Here we use
\[<\xi>^s<\tau + |\xi|^2>^{b+b'}<\xi_1>^{-s}<\xi_2>^{-s + \frac{1}{4} + \epsilon} \le c <\xi>^{-\frac{1}{4}-\epsilon}<\tau_1 - |\xi_1|^2>^{b}\]
to obtain the upper bound
\begin{eqnarray*}
\n{<\xi>^{-\frac{1}{4}-\epsilon} <\tau + |\xi|^2>^{-b}\int d \nu f_1(\xi_1, \tau_1)<\xi_2>^{- \frac{1}{4} - \epsilon}<\tau_2 - |\xi_2|^2>^{-b}f_2(\xi_2, \tau_2)}{L^2_{\xi,\tau}} \\
=c \n{(\F ^{-1}f_1)( J^{s-\frac{1}{4}- \epsilon}u_2)}{\XX{-\frac{1}{4}-\epsilon}{-b}},\hspace{4cm}
\end{eqnarray*}
where $\n{f_1}{L^2_{\xi,\tau}}=\n{\F^{-1}f_1}{L^2_{x,t}}=\n{u_1}{\XX{s}{b}}$. Now we use the dual form of Lemma \ref{l22}, ii), H\"older's inequality and the Lemma itself to obtain
\begin{eqnarray*}
\n{\F ^{-1}f_1 J^{s-\frac{1}{4}- \epsilon}u_2}{\XX{-\frac{1}{4}-\epsilon}{-b}} 
&\le & c\n{\F ^{-1}f_1 J^{s-\frac{1}{4}- \epsilon}u_2}{L_{xt}^{\frac{4}{3}}}\\
&\le & c\n{\F ^{-1}f_1}{L^2_{xt}}\n{J^{s-\frac{1}{4}- \epsilon}u_2}{L^4_{xt}}\\
&\le & c \prod_{i=1}^2 \n{u_i}{\XX{s}{b}}
\end{eqnarray*}
$\hfill \Box$

$Remark:$ In the nonperiodic case we can combine the argument given above with the $L_t^4(L_x^3)$-Strichartz-estimate to obtain the estimate in question whenever $s>-\frac{1}{2}$, $b'< \frac{s}{2}-\frac{1}{4}$, $b>\frac{1}{2}$. (This result was already indicated by Tao, see the remark below Prop. 11.3 in \cite{T00}.) As far as I know, it is still an open question, whether or not the analogue of this Strichartz-estimate, that is
\[\XX{\epsilon}{b} \subset L_t^4({\bf R}, L_x^3({\bf T}^3)),\hspace{1cm}b>\frac{1}{2},\,\,\epsilon>0\]
holds in the periodic case. This, of course, could be used to lower the bound on $s$ in the above theorem down to $-\frac{1}{2}+\epsilon$.

\vspace{0.4cm}

Before we turn to the cubic nonlinearities in the continuous case, let us briefly discuss some counterexamples concerning the periodic case: The examples given by Kenig, Ponce and Vega connected with the onedimensional periodic case (see the proof of Thm 1.10, parts (ii) and (iii) in \cite{KPV96}) show that the estimate
\[\n{u_1 \overline{u}_2}{\XX{s}{b'}}\le c \n{u_1}{\XX{s}{b}}\n{u_2}{\XX{s}{b}}\]
fails for all $s<0$, $b,b' \in {\bf R}$, and that the estimate
\[\n{\overline{u}_1 \overline{u}_2}{\XX{s}{b'}}\le c \n{u_1}{\XX{s}{b}}\n{u_2}{\XX{s}{b}}\]
fails for all $s<- \frac{1}{2}$, if $b-b' \le 1$. From this we can conclude by the method of descent, that these estimates also fail in higher dimensions. So our estimate on $\overline{u}_1 \overline{u}_2$ is sharp (up to the endpoint), while in three dimensions the estimate might be improved (as indicated above), and for $u_1 \overline{u}_2$ no results with $s<0$ can be achieved by the method. For the bilinear form $B(u_1,u_2)=u_1 u_2$ in the two- and threedimensional periodic setting we have the following counterexample exhibiting a significant difference between the periodic and nonperiodic case (cf. the results in \cite{CDKS} and \cite{T00} mentioned in the introduction):

\begin{bsp}\label{ex41}
In the periodic case in space dimension $d\ge 2$ the estimate
\[\n{\prod_{i=1}^2 u_i}{\XX{s}{b'}} \le c \prod_{i=1}^2 \n{u_i}{\XX{s}{b}}\]
fails for all $s<0$, $b, b' \in {\bf R}$.
\end{bsp}

Proof: The above estimate implies
\[\n{<\!\!\tau \!\!+\!\! |\xi|^2\!\!>^{b'} <\!\!\xi\!\!>^{s}\int d \nu  \prod_{i=1}^2 <\!\!\tau_i\!\! + \!|\xi_i|^2\!\!>^{-b} <\!\!\xi_i\!\!>^{-s}f_i(\xi_i,\tau_i)}{L^2_{\xi,\tau}} \le c \prod_{i=1}^2 \n{f_i}{L^2_{\xi,\tau}}.\]
Choosing two orthonormal vectors $e_1$ and $e_2$ in ${\bf R}^d$ and defining for $n \in {\bf N}$
\[f^{(n)}_{1}(\xi,\tau)=\delta_{\xi,ne_1}\chi(\tau + n^2),\hspace{0.3cm}f^{(n)}_2(\xi,\tau)=\delta_{\xi,ne_2}\chi(\tau + n^2),\]
where $\chi$ is the characteristic function of $[-1,1]$, we have $\n{f^{(n)}_i}{L^2_{\xi,\tau}}=c$  and it would follow that
\begin{equation}\label{411}
n^{-2s} \n{<\!\!\tau \!\!+\!\! |\xi|^2\!\!>^{b'} <\!\!\xi\!\!>^{s}\int d \nu \prod_{i=1}^2 f^{(n)}_i(\xi_i,\tau_i)}{L^2_{\xi,\tau}} \le c.
\end{equation}
Now a simple computation shows that
\[\int d \nu \prod_{i=1}^2 f^{(n)}_i(\xi_i,\tau_i) \ge \delta_{\xi,n(e_1 +e_2)}\chi(\tau + 2 n^2),\]
which inserted into (\ref{411}) gives $n^{-s} \le c$. This is a contradiction for all $s<0$.

$\hfill \Box$

\vspace{0.4cm}

The next example shows that our estimate on $\overline{u}_1\overline{u}_2\overline{u}_3$ is essentially sharp:

\begin{bsp}\label{ex42}
In the periodic case in one space dimension the estimate
\[\n{\prod_{i=1}^3 \overline{u}_i}{\XX{s}{b'}} \le c \prod_{i=1}^3\n{u_i}{\XX{s}{b}}\]
fails for all $s< -\frac{1}{3}$, if $b- b' \le 1$.
\end{bsp}

Proof: From the above estimate we obtain
\[\n{<\!\!\tau \!\!+\!\! \xi^2\!\!>^{b'} <\!\!\xi\!\!>^{s}\int d \nu  \prod_{i=1}^3 <\!\!\tau_i\!\! - \!\xi_i^2\!\!>^{-b} <\!\!\xi_i\!\!>^{-s}f_i(\xi_i,\tau_i)}{L^2_{\xi,\tau}} \le c \prod_{i=1}^3 \n{f_i}{L^2_{\xi,\tau}}.\]
Then for $n \in {\bf N}$ we define
\[f^{(n)}_{1,2}(\xi,\tau)=\delta_{\xi,n}\chi(\tau - n^2),\hspace{0.3cm}f^{(n)}_3(\xi,\tau)=\delta_{\xi,-2n}\chi(\tau - 4n^2),\]
with $\chi$ as in the previous example. Again we have $\n{f^{(n)}_i}{L^2_{\xi,\tau}}=c$ and
\begin{equation}\label{412}
n^{-3s} \n{<\!\!\tau \!\!+\!\! \xi^2\!\!>^{b'} <\!\!\xi\!\!>^{s}\int d \nu \prod_{i=1}^3 f^{(n)}_i(\xi_i,\tau_i)}{L^2_{\xi,\tau}} \le c.
\end{equation}
Now it can be easily checked that
\[\int d \nu \prod_{i=1}^3 f^{(n)}_i(\xi_i,\tau_i) \ge \delta_{\xi,0}\chi (\tau - 6n^2).\]
This leads to $n^{-3s+2b'} \le c$ respectively to $\frac{2}{3}b' \le s$.
Consider next the following sequences of functions
\[g^{(n)}_{1}(\xi,\tau)=\delta_{\xi,n}\chi(\tau + 5n^2),\hspace{0.3cm}g^{(n)}_2(\xi,\tau)=\delta_{\xi,n}\chi(\tau - n^2),\hspace{0.3cm}g^{(n)}_3(\xi,\tau)=\delta_{\xi,-2n}\chi(\tau - 4n^2).\]
Arguing as before we are lead to the restriction $-\frac{2}{3}b \le s$. Adding up these two restrictions and taking into account that $b-b' \le 1$ we arrive at $s\ge -\frac{1}{3}$.
$\hfill \Box$

\vspace{0.4cm}

For all the other cubic nonlinearities the corresponding estimates fail for $s<0$, $b, b' \in {\bf R}$, see the examples \ref{ex51} and \ref{ex52} in the next section as well as the remarks below. Next we consider the cubic nonlinearities in the continuous case:

\begin{satz}\label{t111} In the nonperiodic case in one space dimension the estimates
\begin{equation}\label{42}
\n{ \prod_{i=1}^3 \overline{u}_i}{\XX{\sigma}{b'}} \leq c \prod_{i=1}^3 \n{u_i}{\XX{s}{b}}
\end{equation}
and
\begin{equation}\label{43}
\n{\prod_{i=1}^3 u_i}{\XX{\sigma}{b'}} \leq c \prod_{i=1}^3 \n{u_i}{\XX{s}{b}}
\end{equation}
hold, provided $0\ge s>-\frac{5}{12}$, $-\frac{1}{2} < b' < \frac{1}{2}(\frac{1}{4}+ 3s)$, $\sigma <\min{(0,3s-2b')}$ and $b>\frac{1}{2}$.
\end{satz}

Proof: 1. To show (\ref{42}), we write $f_i(\xi, \tau)= <\tau +\xi^2>^{b}<\xi>^s \F u_i(\xi, \tau)$, $1\le i\le 3$. Then we have
\begin{eqnarray*}
&& \n{ \prod_{i=1}^3 \overline{u}_i}{\XX{\sigma}{b'}} =\n{\prod_{i=1}^3 u_i}{\XXm{\sigma}{b'}}\\
& = & \n{<\tau - \xi^2>^{b'} <\xi>^{\sigma}\int d \nu  \prod_{i=1}^3 <\tau_i +\xi_i^2>^{-b} <\xi_i>^{-s}f_i(\xi_i,\tau_i)}{L^2_{\xi,\tau}}.
\end{eqnarray*}
For $0 \le \alpha, \beta, \gamma$ with $\alpha + \beta + \gamma = 2$  we have the inequality
\[<\!\xi_1\!>^{\alpha}<\!\xi_2\!>^{\beta}<\!\xi_3\!>^{\gamma} \le <\!\xi\!>^{2}+ \sum_{i=1}^3 <\!\xi_i\!>^{2}\le c(<\!\tau - \xi^2\!>+\sum_{i=1}^3<\!\tau_i +\xi_i^2\!> \chi_{A_i}),\]
where in $A_i$ the expression $<\tau_i +\xi_i^2>$ is dominant. Hence
\[\n{ \prod_{i=1}^3 \overline{u}_i}{\XX{\sigma}{b'}} \leq c\sum_{k=0}^3 N_k\]
with
\begin{eqnarray*}
N_0 &=& \n{ <\xi>^{\sigma}\int d \nu  \prod_{i=1}^3 <\tau_i +\xi_i^2>^{-b} <\xi_i>^{\frac{2b'}{3}-s}f_i(\xi_i,\tau_i)}{L^2_{\xi,\tau}}\\
&=& c\n{\prod_{i=1}^3 J^\frac{2b'}{3}u_i}{L_t^2(H_x^{\sigma})} \le c\prod_{i=1}^3 \n{J^\frac{2b'}{3}u_i}{\XX{\frac{\sigma}{3}}{b}} \le c\prod_{i=1}^3 \n{u_i}{\XX{s}{b}},
\end{eqnarray*}
where we have used Lemma \ref{l33} and the assumption $\sigma \le 3s-2b' $. Next we estimate $N_1$ by
\begin{eqnarray*}
  \n{<\tau - \xi^2>^{b'} <\xi>^{\sigma}\int d \nu  \prod_{i=1}^3 <\tau_i +\xi_i^2>^{-b} <\xi_i>^{-s}f_i(\xi_i,\tau_i)\chi_{A_1}}{L^2_{\xi,\tau}}\\
\le c \n{\!\!<\!\!\tau \!\!-\!\! \xi^2\!\!>^{-b} <\!\!\xi\!\!>^{\sigma}\int d \nu <\!\!\xi_1\!\!>^{2b'-3s} f_1(\xi_1,\tau_1)\prod_{i=2}^3 <\!\!\tau_i \!\!+\!\!\xi_i^2\!\!>^{-b} f_i(\xi_i,\tau_i)}{L^2_{\xi,\tau}}\\
= c \n{(\Lambda ^b J^{2b'-2s}u_1)(J^su_2)(J^su_3)}{\XXm{\sigma}{-b}},\hspace{3cm}
\end{eqnarray*}
where $\Lambda ^b = \F ^{-1}<\tau + \xi^2>^b \F$. By part iv) of Corollary \ref{k31} this is bounded by
\begin{eqnarray*}
&& c\n{\Lambda ^b J^{2b'-2s}u_1}{L^2_t(H_x^{\sigma})}\n{u_2}{\XX{s}{b}}\n{u_3}{\XX{s}{b}}\\
&=&c \n{u_1}{\XX{2b'-2s + \sigma}{b}}\n{u_2}{\XX{s}{b}}\n{u_3}{\XX{s}{b}} \leq c \prod_{i=1}^3 \n{u_i}{\XX{s}{b}},
\end{eqnarray*}
since $2b'-2s + \sigma \le s$. To estimate $N_k$ for $k=2,3$ one only has to exchange the indices $1$ and $k$. Now (\ref{42}) is shown.

2. Now we prove the second estimate: With $f_i$ as above we have
\[\n{\prod_{i=1}^3 u_i}{\XX{\sigma}{b'}}=c\n{<\!\!\tau \!\!+\!\! \xi^2\!\!>^{b'} <\!\!\xi\!\!>^{\sigma}\int d \nu  \prod_{i=1}^3 <\!\!\tau_i\!\! +\!\!\xi_i^2\!\!>^{-b} <\!\!\xi_i\!\!>^{-s}f_i(\xi_i,\tau_i)}{L^2_{\xi,\tau}}.\]
Here the quantity, which can be controlled by the expressions $<\tau + \xi^2>$, $<\tau_i + \xi_i^2>$, $1\le i \le 3$, is
\[c.q.:= |\xi_1^2 + \xi_2^2+ \xi_3^2 - \xi^2|.\]
So we devide the domain of integration into two parts $A$ and $A^c$, where in $A$ it should hold that
\[\xi_1^2 + \xi_2^2+ \xi_3^2 + \xi^2 \le c\,\,\,\, c.q.\]
Then concerning this region we can argue precisely as in the first part of this proof. For the region $A^c$ we may assume by symmetry that $\xi_1^2 \ge \xi_2^2 \ge \xi_3^2$. Then it is easily checked that in $A^c$ we have
\[1.\,\,\,\xi^2 \ge \frac{1}{2}\xi_1^2 \ge \frac{1}{2}\xi_2^2\,\,\,\hspace{1,5cm}\mbox{and}\hspace{1,5cm}\,\,\,2.\,\,\,\xi_3^2 \le \xi_1^2 \le c (\xi_1 \pm \xi_3)^2.\]
From this it follows
\[\prod_{i=1}^3<\xi_i>^{-s} \le c <\xi>^{-\sigma}<\xi_1 + \xi_3>^{-s_0}<\xi_1 - \xi_3>^{\frac{1}{2}}\]
for $s_0 = \frac{1}{2}+2b' +\epsilon$, so that $-3s \le -\sigma -s_0 + \frac{1}{2} = -\sigma -2b'-\epsilon$ for $\epsilon$ sufficiently small. Hence
\begin{eqnarray*}
&&\n{<\tau + \xi^2>^{b'} <\xi>^{\sigma}\int d \nu  \prod_{i=1}^3 <\tau_i +\xi_i^2>^{-b} <\xi_i>^{-s}f_i(\xi_i,\tau_i)\chi_{A^c}}{L^2_{\xi,\tau}}\\
& \le & c \n{<\!\!\tau\!\! +\!\! \xi^2\!\!>^{b'} \int d \nu <\!\!\xi_1\!\!+\!\!\xi_3\!\!>^{-s_0}<\!\!\xi_1\!\!-\!\!\xi_3\!\!>^{\frac{1}{2}}  \prod_{i=1}^3 <\!\!\tau_i\!\! +\!\!\xi_i^2\!\!>^{-b} f_i(\xi_i,\tau_i)}{L^2_{\xi,\tau}}\\
&=& c \n{(J^su_2) J^{-s_0}J_-^{\frac{1}{2}}(J^su_1,J^su_3)}{\XX{0}{b'}}
\end{eqnarray*}
Using part iii) of Lemma \ref{l23} (observe that $b'<-\frac{1}{4}+\frac{s_0}{2}$) and part i) of Corollary \ref{k23} this can be estimated by
\begin{eqnarray*}
&& c\n{J^su_2}{\XX{0}{b}}\n{J^{-s_0}J_-^{\frac{1}{2}}(J^su_1,J^su_3)}{L^2_t(H_x^{s_0})}\\
& \le & c \n{u_2}{\XX{s}{b}}\n{J_-^{\frac{1}{2}}(J^su_1,J^su_3)}{L^2_{xt}} \le c \prod_{i=1}^3 \n{u_i}{\XX{s}{b}}.
\end{eqnarray*}
$\hfill \Box$

\begin{satz}\label{t1111} In the nonperiodic case in one space dimension the estimate
\begin{equation}\label{11}
\n{ u_1\prod_{i=2}^3 \overline{u}_i}{\XX{s}{b'}} \leq c \prod_{i=1}^3 \n{u_i}{\XX{s}{b}} \nonumber
\end{equation}
holds, provided $-\frac{1}{4} \ge s>-\frac{2}{5}$, $-\frac{1}{2} < b' < \min{(s-\frac{1}{10},-\frac{1}{4}+\frac{s}{2} )}$ and $b>\frac{1}{2}$.
\end{satz}

Proof: We write $f_1(\xi, \tau)= <\tau +\xi^2>^{b}<\xi>^s \F u_1(\xi, \tau)$ and \\ $f_{2,3}(\xi, \tau)= <\tau -\xi^2>^{b}<\xi>^s \F \overline{u}_{2,3}(\xi, \tau)$. Then, using the abbreviations $\sigma_0 = \tau +\xi^2$, $\sigma_1 = \tau_1 +\xi_1^2$ and $\sigma_{2,3} = \tau_{2,3} -\xi_{2,3}^2$, we have
\[\n{ u_1\prod_{i=2}^3 \overline{u}_i}{\XX{s}{b'}}=c\n{<\!\!\sigma_0\!\!>^{b'} <\!\!\xi\!\!>^{s}\int d \nu  \prod_{i=1}^3 <\!\!\sigma_i\!\!>\!\!^{-b} <\!\!\xi_i\!\!>\!\!^{-s}f_i(\xi_i,\tau_i)}{L^2_{\xi,\tau}}.\]
Here the quantity
\[c.q.:= |\xi^2 + \xi_2^2+ \xi_3^2 - \xi_1^2|=2|\xi_2 \xi_3 - \xi (\xi_2+\xi_3)|\]
can be controlled by the expressions $<\sigma_i>$, $0\le i \le 3$. Thus we devide the domain of integration into $A + A^c$, where in $A$ it should hold that $c.q. \ge c <\xi_2><\xi_3>$.
\vspace{0.4cm}

First we consider the region $A^c$. In this region it holds that
\[1. <\xi_2> \le c<\xi>\hspace{1cm} \mbox{or}\hspace{1cm}<\xi_3> \le c<\xi>\]
\[\mbox{and} \hspace{1cm}2.<\xi_{2,3}> \le c<\xi_2 \pm \xi_3>\hspace{0.5cm} \mbox{or}\hspace{0.5cm}<\xi_{2,3}>\le c<\xi \pm \xi_{2,3}>.\]
Writing $A^c=B_1 + B_2$, where in $B_1$ we assume $<\xi_2> \le <\xi_3>$ and in $B_2$, consequently, $<\xi_2> \ge <\xi_3>$, it will be sufficient by symmetry to consider the subregion $B_1$. Now $B_1$ is splitted again into $B_{11}$ and $B_{12}$, where in $B_{11}$ we assume $<\xi_{2,3}> \le c<\xi_2 \pm \xi_3>$ and in $B_{12}$ it should hold that $<\xi_{2,3}>\le c<\xi \pm \xi_{2,3}>$.
\vspace{0.4cm}

Subregion $B_{11}$: Here it holds that \\ $<\xi_1><\xi_2><\xi_3>\le c <\xi><\xi_2-\xi_3><\xi_2+\xi_3>$, giving the upper bound
\begin{eqnarray*}
&&\n{<\sigma_0>^{b'} \int d \nu <\xi_2+\xi_3>^{-s} <\xi_2-\xi_3>^{-s}\prod_{i=1}^3 <\sigma_i>^{-b} f_i(\xi_i,\tau_i)}{L^2_{\xi,\tau}}\\
&=& c \n{(J^s u_1) J^{-s}J_-^{-s}(J^s \overline{u}_2 ,J^s \overline{u}_3) }{\XX{0}{b'}}\\
&\le & c \n{u_1}{\XX{s}{b}}\n{J_-^{-s}(J^s \overline{u}_2 ,J^s \overline{u}_3)}{L^2_{xt}} \le c \prod_{i=1}^3 \n{u_i}{\XX{s}{b}},
\end{eqnarray*}
where we have used part iii) of Lemma \ref{l23} (demanding for $b' < -\frac{1}{4}+ \frac{s}{2}$) and part i) of Corollary \ref{k23}.

\vspace{0.4cm}

Subregion $B_{12}$: Here we have \\ 
$<\xi_1><\xi_2><\xi_3>\le c <\xi><\xi-\xi_3><\xi+\xi_3>$, leading to the upper bound
\begin{eqnarray*}
&&\n{<\!\sigma_0\!>^{b'} \!\int \!d \nu <\!\xi_1+\xi_2+2\xi_3\!>^{-s} <\!\xi_1+\xi_2\!>^{-s}\!\prod_{i=1}^3 \!<\!\sigma_i\!>^{-b} \!f_i(\xi_i,\tau_i)}{L^2_{\xi,\tau}}\\
&=& c \n{J_+^{-s}(J^{-s}((J^s u_1) (J^s \overline{u}_2)) ,J^s \overline{u}_3) }{\XX{0}{b'}}\\
&\le & c \n{u_3}{\XX{s}{b}}\n{J^{-s}((J^s u_1)(J^s \overline{u}_2))}{L^2_{xt}} \le c \prod_{i=1}^3 \n{u_i}{\XX{s}{b}}.
\end{eqnarray*}
Here we have used part ii) of Corollary \ref{k23} (leading again to the restriction \\
$b' < -\frac{1}{4}+ \frac{s}{2}$) and part i) of Lemma \ref{l23}. By this the discussion for the region $A^c$ is completed.

\vspace{0.4cm}

Next we consider the region $A=\sum_{j=0}^3A_j$, where in $A_j$ the expression $<\!\sigma_j\!>$ is assumed to be dominant. By symmetry between the second and third factor (also in the exceptional region $A^c$) it will be sufficient to show the estimate for the subregions $A_0$, $A_1$ and $A_2$. 

\vspace{0.4cm}

Subregion $A_0$: Here we can use $<\xi_2><\xi_3>\le c<\sigma_0>$ to obtain the upper bound
\begin{eqnarray*}
&& \n{ <\!\!\xi\!\!>^{s}\int d \nu <\!\!\sigma_1\!\!>^{-b} <\!\!\xi_1\!\!>^{-s}f_1(\xi_1,\tau_1) \prod_{i=2}^3 <\!\!\sigma_i\!\!>^{-b} <\!\!\xi_i\!\!>^{b'-s}f_i(\xi_i,\tau_i)}{L^2_{\xi,\tau}}\\
&=& c \n{u_1 J^{b'}\overline{u}_2 J^{b'}\overline{u}_3}{L^2_t(H^s_x)}\le c\prod_{i=1}^3 \n{u_i}{\XX{s}{b}}
\end{eqnarray*}
by part ii) of Lemma \ref{l32}, provided $s>-\frac{1}{2}$ (in the last step we have also used $s\ge b'$).

\vspace{0.4cm}

Subregion $A_1$: Here we have $<\xi_2><\xi_3>\le c<\sigma_1>$ and $<\sigma_0> \le <\sigma_1>$. Subdevide $A_1$ again into $A_{11}$ and $A_{12}$ with $<\xi_1>\le c<\xi>$ in $A_{11}$ and, consequently, $<\xi_1>\approx<\xi_2+\xi_3>$ in $A_{12}$. Then for $A_{11}$ we have the upper bound
\begin{eqnarray*}
\n{<\sigma_0>^{-b}\int d \nu f_1(\xi_1,\tau_1) \prod_{i=2}^3 <\sigma_i>^{-b} <\xi_i>^{b'-s}f_i(\xi_i,\tau_i)}{L^2_{\xi,\tau}}\\
=c \n{(\F^{-1}f_1)(J^{b'}\overline{u}_2)(J^{b'}\overline{u}_3)}{\XX{0}{-b}}\le c\n{(\F^{-1}f_1)(J^{b'}\overline{u}_2)(J^{b'}\overline{u}_3)}{L^1_t(L_x^2)}
\end{eqnarray*}
by Sobolev's embedding theorem (plus duality) in the time variable. Now using H\"older's inequality and the $L_t^4(L_x^{\infty})$-Strichartz estimate this can be controlled by
\[\n{\F^{-1}f_1}{L^2_{xt}}\n{J^{b'}u_2}{L_t^4(L_x^{\infty})}\n{J^{b'}u_3}{L_t^4(L_x^{\infty})}\le c\prod_{i=1}^3 \n{u_i}{\XX{s}{b}},\]
provided $b'\le s$.

Now $A_{12}$ is splitted again into $A_{121}$, where we assume $<\xi_2+\xi_3>\le c <\xi_2-\xi_3>$, implying that also $<\xi_1>\le c <\xi_2-\xi_3>$, and $A_{122}$, where $<\xi_2>\approx<\xi_3>$. Consider the subregion $A_{121}$ first: Using $<\xi_1>^{-s}\le c <\xi_2-\xi_3>^{\frac{1}{2}}<\xi_2+\xi_3>^{-s-\frac{1}{2}}$, for this region we obtain the upper bound
\begin{eqnarray*}
\n{\!\!<\!\!\sigma_0\!\!>\!\!\!\!^{-b}\!\!<\!\!\xi\!\!>\!\!\!^{s}
\!\! \int \!\! d \nu f_1(\xi_1,\tau_1)<\!\!\xi_2\!\!-\!\!\xi_3\!\!>\!\!^{\frac{1}{2}}
<\!\!\xi_2\!\!+\!\!\xi_3\!\!>\!\!\!\!^{-s-\frac{1}{2}}
\!\!\prod_{i=2}^3\!\!<\!\!\sigma_i\!\!>\!\!\!\!^{-b}
 \!\!<\!\!\xi_i\!\!>\!\!\!^{b'-s}\!f_i(\xi_i,\tau_i)}{L^2_{\xi,\tau}}&&\\
=  c\n{(\F^{-1}f_1) J^{-s-\frac{1}{2}}J_-^{\frac{1}{2}}(J^{b'}\overline{u}_2,J^{b'}\overline{u}_3)}{\XX{s}{-b}}\hspace{4cm}&&\\
\le  c\n{(\F^{-1}f_1) J^{-s-\frac{1}{2}}J_-^{\frac{1}{2}}(J^{b'}\overline{u}_2,J^{b'}\overline{u}_3)}{L_t^1(L_x^p)}\hspace{1cm}(s-\frac{1}{2}=-\frac{1}{p})\hspace{1.5cm}&&\\
\le  c \n{\F^{-1}f_1}{L^2_{xt}}\n{J^{-s-\frac{1}{2}}J_-^{\frac{1}{2}}(J^{b'}\overline{u}_2,J^{b'}\overline{u}_3)}{L_t^2(L_x^q)}\hspace{1cm}(\frac{1}{p}=\frac{1}{2}+\frac{1}{q})\hspace{1.2cm}&&\\
\le  c \n{u_1}{\XX{s}{b}}\n{J_-^{\frac{1}{2}}(J^{b'}\overline{u}_2,J^{b'}\overline{u}_3)}{L^2_{xt}}\leq c \prod_{i=1}^3 \n{u_i}{\XX{s}{b}}.\hspace{3.5cm}&&
\end{eqnarray*}
Next we consider the subregion $A_{122}$, where $<\xi_2>\approx<\xi_3> \ge c <\xi_1>$. Here we get the upper bound
\begin{eqnarray*}
&& \n{\!\!<\!\!\sigma_0\!\!>\!\!\!\!^{-b}\!\!<\!\!\xi\!\!>\!\!\!^{s}
\!\! \int \!\! d \nu f_1(\xi_1,\tau_1)
<\!\!\xi_1\!\!>\!\!\!\!^{-s-\frac{1}{6}}\!\!\prod_{i=2}^3 \!\!<\!\!\sigma_i\!\!>\!\!\!\!^{-b}\!\!<\!\!\xi_i\!\!>\!\!\!^{b'-s+\frac{1}{12}} \!f_i(\xi_i,\tau_i)}{L^2_{\xi,\tau}}\\
& = & c \n{(\Lambda ^b J ^{-\frac{1}{6}}u_1)(J^{b'+\frac{1}{12}}\overline{u}_2)(J^{b'+\frac{1}{12}}\overline{u}_3)}{\XX{s}{-b}},\hspace{1cm}(\Lambda ^b = \F^{-1}<\tau + \xi^2>^b\F)\\
& \le & c \n{\Lambda ^b u_1}{L^2_t(H_x^{-\frac{1}{4}-\frac{1}{6}})}\n{J^{b'+\frac{1}{12}}u_2}{\XX{0}{b}}\n{J^{b'+\frac{1}{12}}u_3}{\XX{0}{b}},
\end{eqnarray*}
where we have used $s \le -\frac{1}{4}$ and part iv) of Corollary \ref{k31}. The latter is bounded by $c \prod_{i=1}^3 \n{u_i}{\XX{s}{b}}$, provided $s \ge -\frac{5}{12}$ and $s \ge b' +\frac{1}{12}$. Thus the discussion for the region $A_1$ is complete.

\vspace{0.4cm}

Subregion $A_2$: First we write $A_2=A_{21}+A_{22}$, where in $A_{21}$ it should hold that $<\xi_1> \le c <\xi>$. Then this subregion can be treated precisely as the subregion $A_{11}$, leading to the bound $s > - \frac{1}{2}$. For the remaining subregion $A_{22}$ it holds that
\[<\xi_2><\xi_3>\le c <\sigma_2>\hspace{1cm}\mbox{and}\hspace{1cm} <\xi_1>\le c <\xi_2 + \xi_3>.\]
Now $A_{22}$ is splitted again into $A_{221}$, where we assume $<\xi_1>\,\,\,\le c \,\,\,<\xi_2>$, and into $A_{222}$, where we then have $<\xi_2> \,\,\,<< \,\,\,<\xi_1>$. The upper bound for $A_{221}$ is
\begin{eqnarray*}
&& \n{\!\!<\!\!\sigma_0\!\!>\!\!\!\!^{-b}\!\!<\!\!\xi\!\!>\!\!\!^{s}
\!\! \int \!\! d \nu f_2(\xi_2,\tau_2)
<\!\!\xi_2\!\!>\!\!\!\!^{-s}\!\prod_{i\ne 2} \!\!<\!\!\sigma_i\!\!>\!\!\!\!^{-b}\!\!<\!\!\xi_i\!\!>\!\!\!^{b'-s} \!f_i(\xi_i,\tau_i)}{L^2_{\xi,\tau}}\\
& \le & c \n{(\Lambda ^b_-\overline{u}_2)(J^{b'}u_1)(J^{b'}\overline{u}_3)}{\XX{s}{-b}} \hspace{1cm}(\Lambda ^b_- = \F ^{-1} <\tau - \xi ^2>^b \F)\\
& \le & c \n{\Lambda ^b_-\overline{u}_2}{L_t^2(H_x^s)}\n{u_1}{\XX{b'}{b}}\n{u_3}{\XX{b'}{b}}\le c \prod_{i=1}^3 \n{u_i}{\XX{s}{b}}.
\end{eqnarray*}
Here we have used part i) of Lemma \ref{l32} (dualized version) and the assumption $s \ge b'$.

For the subregion $A_{222}$ the argument is a bit more complicated and it is here, where the strongest restrictions on $s$ occur: Subdevide $A_{222}$ again into $A_{2221}$ and $A_{2222}$ with $<\xi_2>^2 \le <\xi_1>$ in $A_{2221}$. Then in $A_{2221}$ it holds that
\[(<\xi_1><\xi_2><\xi_3>)^{\frac{2}{5}} \le c <\xi_1> \le c <\xi_3>\le c<\xi_2 \pm \xi_3>,\]
hence, for $\epsilon = 1 + \frac{5}{2}s\,\,\,(>0)$,
\[\prod_{i=1}^3<\xi_i>^{-s} \le c <\xi_2-\xi_3>^{\frac{1}{2}}<\xi_2+\xi_3>^{\frac{1}{2}-\epsilon}.\]
Then, throwing away the $<\xi>^s$-factor, we obtain the upper bound
\begin{eqnarray*}
&& \n{\!\!<\!\!\sigma_0\!\!>\!\!\!^{b'}\!\!<\!\!\xi_2\!\!-\!\!\xi_3\!\!>\!\!^{\frac{1}{2}}<\!\!\xi_2\!\!+\!\!\xi_3\!\!>^{\frac{1}{2}-\epsilon}\prod_{i=1}^3\!\!<\!\!\sigma_i\!\!>\!\!\!\!^{-b}f_i(\xi_i,\tau_i)}{L^2_{\xi,\tau}}\\
& = & c \n{(J^su_1)J^{\frac{1}{2}-\epsilon}J_-^{\frac{1}{2}}(J^{s}\overline{u}_2,J^{s}\overline{u}_3)}{\XX{0}{b'}}\\
& \le & c \n{u_1}{\XX{s}{b}} \n{J_-^{\frac{1}{2}}(J^{s}\overline{u}_2,J^{s}\overline{u}_3)}{L^2_{xt}} \le c \prod_{i=1}^3 \n{u_i}{\XX{s}{b}},
\end{eqnarray*}
by Lemma \ref{l23}, part iii), and Corollary \ref{k23}, part i) (and the remark below), leading to the restriction $b' < \frac{5}{4}s$, which - in the allowed range for $s$ - is in fact weaker than $b' < s - \frac{1}{10}$.
Finally we consider the subregion $A_{2222}$, where we have $<\xi_1>^{\frac{1}{2}}\le <\xi_2>\,\,\,<<\,\,\,<\xi_1>\approx<\xi_3>$, implying that
\[<\xi_1>^{\frac{3}{20}} \le c (<\xi_2><\xi_3>)^{\frac{1}{10}}.\]
This gives the upper bound
\begin{eqnarray*}
 \n{\!\!<\!\!\sigma_0\!\!>\!\!\!\!^{-b}\!\!<\!\!\xi\!\!>\!\!\!^{s}
\!\! \int \!\! d \nu <\!\!\xi_1\!\!>\!\!\!\!^{-s-\frac{3}{20}}\!<\!\!\sigma_1\!\!>\!\!\!\!^{-b}f_1(\xi_1,\tau_1)\!\!\prod_{i=2}^3 \!\!<\!\!\xi_i\!\!>\!\!\!^{b'-s+\frac{1}{10}}f_i(\xi_i,\tau_i)\!\!<\!\!\sigma_3\!\!>\!\!\!\!^{-b}}{L^2_{\xi,\tau}}\\
\le  c \n{(J^{-\frac{3}{20}}u_1)(\Lambda ^b_- J^{b'+\frac{1}{10}}\overline{u}_2)(J^{b'+\frac{1}{10}}\overline{u}_3)}{\XX{s}{b}}.\hspace{3cm}
\end{eqnarray*}
Now using $s \le -\frac{1}{4}$ again and part ii) of Corollary \ref{k31} this can be estimated by
\[c\n{u_1}{\XX{-\frac{3}{20}-\frac{1}{4}}{b}}\n{\Lambda ^b_- J^{b'+\frac{1}{10}}\overline{u}_2}{L^2_{xt}}\n{u_3}{\XX{b'+\frac{1}{10}}{b}}\le c \prod_{i=1}^3 \n{u_i}{\XX{s}{b}},\]
since $s>-\frac{2}{5}$ and $s>b'+\frac{1}{10}$ as assumed.
$\hfill \Box$

\vspace{0.4cm}

$Remark$: The estimate (\ref{11}) also holds under the assumption $s \ge -\frac{1}{4}$, $b'<-\frac{3}{8}$ and $b>\frac{1}{2}$. For $s = -\frac{1}{4}$ this is contained in the above theorem, and for $s > -\frac{1}{4}$ this follows from $<\xi > \le c \prod_{i=1}^3<\xi_i >$.

\section{Estimates on quartic nonlinearities}

\begin{satz}\label{t2} Let $n=1$. Assume $ 0 \geq s > - \frac{1}{6}$ and $ - \frac{1}{2} < b' < \frac{3s}{2} - \frac{1}{4}$. Then in the periodic and nonperiodic case for all $b>\frac{1}{2}$ the estimate 
\[\n{ \prod_{i=1}^4\overline{u}_i}{\XX{s}{b'}} \leq c \prod_{i=1}^4 \n{u_i}{\XX{s}{b}}\]
holds true.
\end{satz}

Proof: Again we write $f_i(\xi, \tau)= <\tau - \xi^2>^{b}<\xi>^s \F \overline{u}_i(\xi, \tau)$, so that 
\[\n{ \prod_{i=1}^4 \overline{u}_i}{\XX{s}{b'}}=c \n{<\!\!\tau \!\!+\!\! \xi^2\!\!>^{b'} <\!\!\xi\!\!>^{s}\int d \nu  \prod_{i=1}^4 <\!\!\tau_i\!\! - \!\!\xi_i^2\!\!>^{-b} <\!\!\xi_i\!\!>^{-s}f_i(\xi_i,\tau_i)}{L^2_{\xi,\tau}}.\]
Now we can use the inequality (\ref{41}) with $m=4$ and the assumption $b'<\frac{3 s}{2} - \frac{1}{4}$ to obtain
\[<\xi>^{s + \frac{1}{2}-\epsilon}\prod_{i=1}^4<\xi_i>^{-s+\epsilon} \leq c (<\tau + \xi^2 >^{-b'} + \sum_{i=1}^4 <\tau_i - \xi_i^2>^{-b'} \chi_{A_i})\]
for some $\epsilon>0$. (Again in $A_i$ we assume $<\tau_i - \xi_i^2> \,\,\ge \,\, <\tau + \xi^2 >$.) From this it follows, that
\[\n{ \prod_{i=1}^4 \overline{u}_i}{\XX{s}{b'}} \le c \sum_{j=0}^4 \n{I_j}{L^2_{\xi,\tau}},\]
with
\[I_0(\xi,\tau)=<\xi>^{-\frac{1}{2}+\epsilon}\int d \nu  \prod_{i=1}^4<\tau_i - \xi_i^2>^{-b} <\xi_i>^{-\epsilon}f_i(\xi_i,\tau_i)\]
and, for $1 \le j \le m$,
\begin{eqnarray*}
I_j(\xi,\tau)\!\!\!\! & = & \!\!\!\!<\!\!\xi\!\!>^{-\frac{1}{2}+\epsilon}<\!\!\tau \!\!+\!\! \xi^2 \!\!>^{b'}\!\!\int\!\! d \nu <\!\!\tau_j \!\!- \!\!\xi_j^2 \!\! >^{-b'} \prod_{i=1}^4 <\!\!\tau_i\!\! -\!\! \xi_i^2\!\!>^{-b} <\!\!\xi_i\!\!>^{-\epsilon}f_i(\xi_i,\tau_i)\chi_{A_i} \\
\!\!\!\!& \le &\!\!\!\! <\!\!\xi\!\!>^{-\frac{1}{2}+\epsilon}<\!\!\tau + \xi^2 \!\!>^{-b}\!\!\int \!\! d \nu <\!\!\tau_j - \xi_j^2\!\! >^{b} \prod_{i=1}^4 <\!\!\tau_i - \xi_i^2\!\!>^{-b} <\!\!\xi_i\!\!>^{-\epsilon} f_i(\xi_i,\tau_i).
\end{eqnarray*}

Next we estimate $I_0$ using first Sobolev's embedding theorem, then H\"older's inequality, again Sobolev and finally Corollary \ref{k21}. Here $\epsilon'$, $\epsilon''$ denote suitable small, positive numbers.

\begin{eqnarray*}
\n{I_0}{L^2_{\xi,\tau}}
&=&\n{\prod_{i=1}^4 J^{s-\epsilon} \overline{u}_i}{L^2_t (H_x^{-\frac{1}{2}+\epsilon})}
 \le c\n{\prod_{i=1}^4 J^{s-\epsilon} \overline{u}_i}{L^2_t (L_x^{1+\epsilon'})}\\
& \le &c\prod_{i=1}^4\n{ J^{s-\epsilon} \overline{u}_i}{L_t^{8}(L_x^{4+4\epsilon'})}
 \le c\prod_{i=1}^4\n{ J^{s-\epsilon''} \overline{u}_i}{L_t^{8}(L_x^{4})}\\
& \le &c\prod_{i=1}^m\n{ J^{s}\overline{u}_i}{\XXm{0}{b}}=c\prod_{i=1}^m\n{\overline{u}_i}{\XXm{s}{b}}
\end{eqnarray*}
To estimate $I_j$, $1\le j \le m$, we use Sobolev (in both variables) plus duality, H\"older, again Sobolev (in the space variable) and Lemma \ref{l21}. Again we need suitable small, positive numbers $\epsilon'$, $\epsilon''$ and $\epsilon'''$.

\begin{eqnarray*}
\n{I_j}{L^2_{\xi,\tau}} 
& \le & c\n{(\prod_{\stackrel{i=1}{i \ne j}}^4 J^{s-\epsilon} \overline{u}_i)(J^{-\epsilon}\F^{-1}f_j)}{\XX{-\frac{1}{2}+\epsilon}{-b}}\\
& \le & c\n{(\prod_{\stackrel{i=1}{i \ne j}}^4 J^{s-\epsilon} \overline{u}_i)(J^{-\epsilon}\F^{-1}f_j)}{L^{1}_t(L_x^{1+\epsilon'})}\\
& \le & c\n{J^{-\epsilon}\F^{-1}f_j}{L^2_{x,t}}\prod_{\stackrel{i=1}{i \ne j}}^4\n{ J^{s-\epsilon} \overline{u}_i}{L^{6}_t(L_x^{6+\epsilon''})}\\
& \le & c\n{J^{-\epsilon}\F^{-1}f_j}{L^2_{x,t}}\prod_{\stackrel{i=1}{i \ne j}}^4\n{ J^{s-\epsilon'''} \overline{u}_i}{L^{6}_{xt}}\\
& \le & c\n{f_j}{L^2_{\xi,\tau}}\prod_{\stackrel{i=1}{i \ne j}}^4\n{ J^{s} \overline{u}_i}{\XXm{0}{b}} = c \prod_{i=1}^4\n{\overline{u}_i}{\XXm{s}{b}}
\end{eqnarray*}
$\hfill \Box$

\vspace{0.4cm}

In the periodic case the following examples show, that for all the other quartic nonlinearities ($u^4,\,\,u^3\overline{u},\,\,...,\,\,u\overline{u}^3$) the corresponding estimates fail for all $s<0$. The argument is essentially that given in the proof of Thm 1.10 in \cite{KPV96}.

\begin{bsp}\label{ex51}
In the periodic case in one space dimension the estimate
\[\n{\prod_{i=1}^4 u_i}{\XX{s}{b'}} \le c \prod_{i=1}^4 \n{u_i}{\XX{s}{b}}\]
fails for all $s<0$, $b, b' \in {\bf R}$.
\end{bsp}

Proof: The above estimate implies
\[\n{<\!\!\tau \!\!+\!\! \xi^2\!\!>^{b'} <\!\!\xi\!\!>^{s}\int d \nu  \prod_{i=1}^4 <\!\!\tau_i\!\! + \!\!\xi_i^2\!\!>^{-b} <\!\!\xi_i\!\!>^{-s}f_i(\xi_i,\tau_i)}{L^2_{\xi,\tau}} \le c \prod_{i=1}^4 \n{f_i}{L^2_{\xi,\tau}}.\]
Defining for $n \in {\bf N}$
\[f^{(n)}_{1,2}(\xi,\tau)=\delta_{\xi,2n}\chi(\tau + \xi^2),\hspace{0.3cm}f^{(n)}_3(\xi,\tau)=\delta_{\xi,-n}\chi(\tau + \xi^2),\hspace{0.3cm}f^{(n)}_4(\xi,\tau)=\delta_{\xi,0}\chi(\tau + \xi^2),\]
where $\chi$ is the characteristic function of $[-1,1]$, we have $\n{f^{(n)}_i}{L^2_{\xi,\tau}}=c$  and it would follow that
\begin{equation}\label{51}
n^{-3s} \n{<\!\!\tau \!\!+\!\! \xi^2\!\!>^{b'} <\!\!\xi\!\!>^{s}\int d \nu \prod_{i=1}^4 f^{(n)}_i(\xi_i,\tau_i)}{L^2_{\xi,\tau}} \le c.
\end{equation}
Now a simple computation shows that
\[\int d \nu \prod_{i=1}^4 f^{(n)}_i(\xi_i,\tau_i) \ge \delta_{\xi,3n}\chi(\tau + \xi^2).\]
Inserting this into (\ref{51}) we obtain $n^{-2s} \le c$, which is a contradiction for any $s<0$.

$\hfill \Box$

\vspace{0.4cm}

$Remark:$ Using only the sequences $f^{(n)}_{1,2,3}$ from the above proof, the same calculation shows that in the periodic case the estimate
\[\n{\prod_{i=1}^3 u_i}{\XX{s}{b'}} \le c \prod_{i=1}^3 \n{u_i}{\XX{s}{b}}\]
fails for all $s<0$, $b, b' \in {\bf R}$.

\vspace{0.4cm}

\begin{bsp}\label{ex52}
In the periodic case in one space dimension the estimates
\[\n{ u_1 \overline{u}_2 \tilde{u}_3 \tilde{u}_4}{\XX{s}{b'}} \le c \prod_{i=1}^4 \n{u_i}{\XX{s}{b}},\]
where $\tilde{u}=u$ or $\tilde{u}=\overline{u}$, fail for all $s<0$, $b, b' \in {\bf R}$.
\end{bsp}

Proof: We define
\begin{eqnarray*}
f^{(n)}_{1}(\xi,\tau)=\delta_{\xi,n}\chi(\tau + \xi^2)&,&f^{(n)}_{2}(\xi,\tau)=\delta_{\xi,-n}\chi(\tau - \xi^2)\\
f^{(n)}_{3,4}(\xi,\tau)=\delta_{\xi,0}\chi(\tau \pm \xi^2)&&(+ \,\,\,\,\mbox{for}\,\,\,\,\tilde{u}_{3,4}=u_{3,4},\,\,\,- \,\,\,\,\mbox{for}\,\,\,\,\tilde{u}_{3,4}=\overline{u}_{3,4}).
\end{eqnarray*}
Then the above estimate would imply
\begin{equation}\label{52}
n^{-2s} \n{<\!\!\tau \!\!+\!\! \xi^2\!\!>^{b'} <\!\!\xi\!\!>^{s}\int d \nu \prod_{i=1}^4 f^{(n)}_i(\xi_i,\tau_i)}{L^2_{\xi,\tau}} \le c.
\end{equation}
Now
\[\int d \nu \prod_{i=1}^4 f^{(n)}_i(\xi_i,\tau_i) \ge \delta_{\xi,0}\chi(\tau ),\]
which inserted into (\ref{52}) again leads to $n^{-2s} \le c$.
$\hfill \Box$

\vspace{0.4cm}

$Remark:$ Using only the sequences $f^{(n)}_{1,2,3}$ from the above proof, we see that in the periodic case the estimates
\[\n{u_1 \overline{u}_2 \tilde{u}_3}{\XX{s}{b'}} \le c \prod_{i=1}^3 \n{u_i}{\XX{s}{b}}\]
fail for all $s<0$, $b, b' \in {\bf R}$.

\vspace{0.4cm}

Now we turn to discuss the continuous case, where we can use the bi- and trilinear inequalities of section 2.2 respectively 3 in order to prove the relevant estimates for some $s<0$. We start with the following

\begin{prop}\label{p51}
Let $0\ge s > -\frac{1}{8}$, $-\frac{1}{2}<b'<-\frac{1}{4}+2s$. Then in the continuous case in one space dimension for any $b>\frac{1}{2}$ the estimate
\[\n{u_1u_2\overline{u}_3\overline{u}_4}{\XX{s}{b'}}\leq c \prod_{i=1}^4 \n{u_i}{\XX{s}{b}}\]
holds true.
\end{prop}

Proof: Apply part iii) of Lemma \ref{l23} to obtain
\[\n{u_1u_2\overline{u}_3\overline{u}_4}{\XX{s}{b'}}\leq c \n{u_1}{\XX{s}{b}}\n{u_2\overline{u}_3\overline{u}_4}{L_t^2(H^{\sigma - s})},\]
provided that $s\le 0$, $-\frac{1}{2}\le \sigma \le 0$, $b'<-\frac{1}{4}+\frac{\sigma}{2}$. This is fulfilled for $\sigma = 4s$ and the second factor is equal to
\[\n{u_2\overline{u}_3\overline{u}_4}{L_t^2(H^{3 s})} \le c \prod_{i=2}^4 \n{u_i}{\XX{s}{b}}\]
by Lemma \ref{l32} and the remark below.
$\hfill \Box$

\vspace{0.4cm}

To show that this proposition is essentially (up to the endpoint) sharp, we present the following counterexample (cf. Thm 1.4 in \cite{KPV96}):

\begin{bsp}\label{ex53}
In the nonperiodic case in one space dimension the estimate
\[\n{ u_1 u_2 \overline{u}_3 \overline{u}_4}{\XX{s}{b'}} \le c \prod_{i=1}^4 \n{u_i}{\XX{s}{b}}\]
fails for all $s<-\frac{1}{8}$, $b, b' \in {\bf R}$.
\end{bsp}

Proof: The above estimate implies
\[\n{<\!\!\tau \!\!+\!\! \xi^2\!\!>^{b'} <\!\!\xi\!\!>^{s}\int d \nu  \prod_{i=1}^4 <\!\!\sigma_i\!\! >^{-b} <\!\!\xi_i\!\!>^{-s}f_i(\xi_i,\tau_i)}{L^2_{\xi,\tau}} \le c \prod_{i=1}^4 \n{f_i}{L^2_{\xi,\tau}},\]
where $<\sigma_{1,2}>=<\tau_{1,2}+\xi^2_{1,2}>$ and $<\sigma_{3,4}>=<\tau_{3,4}-\xi^2_{3,4}>$. Choosing
\[f^{(n)}_{1,2}(\xi,\tau)=\chi(\xi - n)\chi(\tau + \xi^2),\,\,\,\,\,\,\,\,f^{(n)}_{3,4}(\xi,\tau)=\chi(\xi + n)\chi(\tau - \xi^2)\]
we arrive at
\begin{equation}\label{53}
n^{-4s} \n{<\!\!\tau \!\!+\!\! \xi^2\!\!>^{b'} <\!\!\xi\!\!>^{s}\int d \nu \prod_{i=1}^4 f^{(n)}_i(\xi_i,\tau_i)}{L^2_{\xi,\tau}} \le c.
\end{equation}
Now an elementary computation gives
\[\int d \nu \prod_{i=1}^4 f^{(n)}_i(\xi_i,\tau_i) \ge c\chi_c(2n\xi )\chi_c(\tau ),\]
where $\chi_c$ is the characteristic function of $[-c,c]$. Inserting this into (\ref{53}) we get $n^{-4s-\frac{1}{2}} \le c$, which is a contradiction for any $s<-\frac{1}{8}$.
$\hfill \Box$

Finally we consider the remaining nonlinearities $u^4$, $u^3\overline{u}$ and $u\overline{u}^3$, for which we can lower the bound on $s$ down to $-\frac{1}{6}+ \epsilon$ :

\begin{satz}\label{t22} Let $n=1$. Assume $ 0 \geq s > - \frac{1}{6}$, $ - \frac{1}{2} < b' < \frac{3s}{2} - \frac{1}{4}$ and $b>\frac{1}{2}$. Then in the nonperiodic case the estimates 
\[\n{N(u_1,u_2,u_3,u_4) }{\XX{s}{b'}} \leq c \prod_{i=1}^4 \n{u_i}{\XX{s}{b}}\]
hold true for $N(u_1,u_2,u_3,u_4)=\prod_{i=1}^4u_i$, $=(\prod_{i=1}^3u_i)\overline{u}_4$ or $=(\prod_{i=1}^3\overline{u}_i)u_4$.
\end{satz}

Proof: 1. We begin with the nonlinearity $N(u_1,u_2,u_3,u_4)=\prod_{i=1}^4u_i$. Writing $f_i(\xi, \tau)= <\tau + \xi^2>^{b}<\xi>^s \F u_i(\xi, \tau)$ we have
\[\n{ \prod_{i=1}^4 u_i}{\XX{s}{b'}}=c \n{<\!\!\tau \!\!+\!\! \xi^2\!\!>^{b'} <\!\!\xi\!\!>^{s}\int d \nu  \prod_{i=1}^4 <\!\!\tau_i\!\! +\!\!\xi_i^2\!\!>^{-b} <\!\!\xi_i\!\!>^{-s}f_i(\xi_i,\tau_i)}{L^2_{\xi,\tau}}.\]
The quantity controlled by the expression $<\tau + \xi^2>$, $<\tau_i+\xi_i^2>$, $1\le i\le 4$ is $|\sum_{i=1}^4 \xi_i^2 - \xi^2|$. We devide the domain of integration into $A$ and $A^c$, where in $A$ we assume $\xi^2 \le \frac{\xi_1^2}{2}$ and thus
\[|\sum_{i=1}^4 \xi_i^2 - \xi^2| \ge c (\sum_{i=1}^4 \xi_i^2 + \xi^2).\]
So concerning this region we may refer to the proof of Theorem \ref{t2}. For the region $A^c$, where $\xi_1^2 \le 2 \xi^2$ we have the upper bound
\[c\n{(J^su_1)\prod_{i=2}^4u_i}{\XX{0}{b'}} \le c \n{u_1}{\XX{s}{b}}\n{\prod_{i=2}^4u_i}{L_t^2(H_x^{3s})}\le c \prod_{i=1}^4 \n{u_i}{\XX{s}{b}}\]
by Lemma \ref{l23}, part iii), which requires $b' < \frac{3s}{2} - \frac{1}{4}$, $s\ge - \frac{1}{6}$, and by Lemma \ref{l33} (and the remark below), which demands $s> - \frac{1}{6}$.

\vspace{0.4cm}

2. Next we consider $N(u_1,u_2,u_3,u_4)=(\prod_{i=1}^3u_i)\overline{u}_4$. For $1\le i \le 3$ we choose the $f_i$ as in the first part of this proof and $f_4(\xi, \tau)= <\tau - \xi^2>^{b}<\xi>^s \F \overline{u}_4(\xi, \tau)$. Then the left hand side of the claimed estimate is equal to
\[c \n{\!\!<\!\!\tau \!\!+\!\! \xi^2\!\!>^{b'} <\!\!\xi\!\!>^{s}\!\!\int\!\! d \nu  \prod_{i=1}^3\!\! <\!\!\tau_i\!\! +\!\!\xi_i^2\!\!>\!\!\!\!^{-b} \!\!<\!\!\xi_i\!\!>\!\!\!\!^{-s}f_i(\xi_i,\tau_i)<\!\!\tau_4\!\! -\!\!\xi_4^2\!\!>\!\!\!\!^{-b} \!\!<\!\!\xi_4\!\!>\!\!\!\!^{-s}f_4(\xi_4,\tau_4)}{L^2_{\xi,\tau}}.\]
Now the quantity controlled by the expressions $<\tau + \xi^2>$, $<\tau_i + \xi_i^2>$, $1\le i \le 3$, and $<\tau_4 - \xi_4^2>$ is
\[c.q.:=|\xi_1^2 + \xi_2^2 + \xi_3^2 - \xi_4^2 - \xi^2 |.\]
We devide the domain of integration into the regions $A$, $B$ and $C =(A+B)^c$, where in $A$ it should hold that
\[c.q.\,\,\,\ge c (\sum_{i=1}^4 \xi_i^2 + \xi^2).\]
Again, concerning this region we may refer to the proof of Thm. \ref{t2}.

Next we write $B = \bigcup_{i=1}^3 B_i$, where in $B_i$ we assume $\xi_i^2 \le c \xi^2$ for some large constant $c$. By symmetry it is sufficient to consider the subregion $B_1$, where we obtain the upper bound
\[c \n{(J^su_1)u_2 u_3 \overline{u}_4}{\XX{0}{b'}} \le c \n{u_1}{\XX{s}{b}}\n{u_2 u_3 \overline{u}_4}{L^2_t(H_x^{3s})}\le c \prod_{i=1}^4 \n{u_i}{\XX{s}{b}}\]
by Lemma \ref{l23}, part iii), demanding for $b' < \frac{3s}{2} - \frac{1}{4}$, $3s \ge -\frac{1}{2}$ and Lemma \ref{l32} (resp. the remark below), where $s>-\frac{1}{6}$ is required.

Considering the region $C$ we may assume by symmetry between the first three factors that $\xi_1^2\ge \xi_2^2 \ge \xi_3^2$. Then for this region it is easily checked that

\vspace{0.4cm}

1. $\xi^2 << \xi_3^2$,
\hspace{0.4cm}
2. $\xi_4^2 \ge \frac{3}{2} \xi_2^2$,  hence  $\xi_4^2 \le c (\xi_4 + \xi_2)^2$,  and
\hspace{0.4cm}3. $\xi_1^2 \le c (\xi_1 - \xi_3)^2$.

\vspace{0.4cm}

This implies
\begin{itemize}
\item[1.] $<\xi>^{-2s}<\xi_4>^{-s} \le c <\xi_4 + \xi_2>^{-3s}$\hspace{0.5cm}and
\item[2.] $<\xi>^{\frac{1}{2}+3s}<\xi_1>^{-s}<\xi_2>^{-s}<\xi_3>^{-s}\le c <\xi_1 - \xi_3>^{\frac{1}{2}}$,
\end{itemize}
leading to the upper bound
\begin{eqnarray*}
&& \n{J_-^{\frac{1}{2}}(J^su_1,J^su_3) J^{-3s}(J^su_2J^s\overline{u}_4)}{\XX{-\frac{1}{2}}{b'}}\\
& \le & c \n{J_-^{\frac{1}{2}}(J^su_1,J^su_3) J^{-3s}(J^su_2J^s\overline{u}_4)}{L_t^p(L_x^{1+})}\hspace{2cm}(b'-\frac{1}{2}=-\frac{1}{p})\\
& \le & c \n{J_-^{\frac{1}{2}}(J^su_1,J^su_3)}{L^2_{xt}}\n{J^{-3s}(J^su_2J^s\overline{u}_4)}{L_t^q(L_x^{2+})}\hspace{1cm}(\frac{1}{q}=\frac{1}{p}-\frac{1}{2}=-b').
\end{eqnarray*}
Using Corollary \ref{k23} the first factor can be estimated by
\[c \n{u_1}{\XX{s}{b}}\n{u_3}{\XX{s}{b}},\]
while for the second we can use Sobolev's embedding Theorem and part ii) of Lemma \ref{l23} to obtain the bound
\[c \n{J^su_2J^s\overline{u}_4}{L_t^q(H_x^{-3s+})} \le c \n{u_2}{\XX{s}{b}}\n{u_4}{\XX{s}{b}}.\]
Here the restriction $b' < \frac{3s}{2} - \frac{1}{4}$ is required again.

\vspace{0.4cm}

3. Finally we consider $N(u_1,u_2,u_3,u_4)=(\prod_{i=1}^3\overline{u}_i)u_4$. With $f_i(\xi, \tau)= <\tau - \xi^2>^{b}<\xi>^s \F \overline{u}_i(\xi, \tau)$, $1\le i \le 3$ and $f_4(\xi, \tau)= <\tau + \xi^2>^{b}<\xi>^s \F u_4(\xi, \tau)$
the norm on the left hand side is equal to
\[c \n{\!\!<\!\!\tau \!\!+\!\! \xi^2\!\!>^{b'} <\!\!\xi\!\!>^{s}\!\!\int\!\! d \nu  \prod_{i=1}^3\!\! <\!\!\tau_i\!\! -\!\!\xi_i^2\!\!>\!\!\!\!^{-b} \!\!<\!\!\xi_i\!\!>\!\!\!\!^{-s}f_i(\xi_i,\tau_i)<\!\!\tau_4\!\! +\!\!\xi_4^2\!\!>\!\!\!\!^{-b} \!\!<\!\!\xi_4\!\!>\!\!\!\!^{-s}f_4(\xi_4,\tau_4)}{L^2_{\xi,\tau}}.\]
The controlled quantity here is
\[c.q.:=|\xi_1^2 + \xi_2^2 + \xi_3^2 - \xi_4^2 + \xi^2 |.\]
Devide the area of integration into $A$, $B$ and $C =(A+B)^c$, where in $A$ we assume again
\[c.q.\,\,\,\ge c (\sum_{i=1}^4 \xi_i^2 + \xi^2)\]
in order to refer to the proof of Thm. \ref{t2}. In $B$ we assume $\xi_4^2 \le c \xi^2$, so that for this region we have the bound
\[c \n{\overline{u}_1 \overline{u}_2 \overline{u}_3(J^su_4)}{\XX{0}{b'}} \le c \n{u_4}{\XX{s}{b}}\n{u_2 u_3 u_4}{L^2_t(H_x^{3s})}\le c \prod_{i=1}^4 \n{u_i}{\XX{s}{b}}\]
by Lemma \ref{l23}, part iii), and Lemma \ref{l33} and the remark below. Here $b' < \frac{3s}{2} - \frac{1}{4}$ and $s>-\frac{1}{6}$ is required.

For the region $C$ we shall assume again that $\xi_1^2\ge \xi_2^2 \ge \xi_3^2$. Then it is easily checked that in $C$

\vspace{0.4cm}

1. $\xi^2 << \xi_4^2$,
\hspace{0.4cm}
2. $\xi_4^2 \ge \frac{3}{2} \xi_2^2$,  hence  $\xi_4^2 \le c (\xi_4 + \xi_2)^2$,  and
\hspace{0.4cm}3. $\xi_1^2 \le c (\xi_1 - \xi_3)^2$.

\vspace{0.4cm}

Then for $C$ we have the estimate
\begin{eqnarray*}
&& \n{J_-^{\frac{1}{2}}(J^s\overline{u}_1,J^s\overline{u}_3) J^{-3s}(J^su_4J^s\overline{u}_2)}{\XX{-\frac{1}{2}}{b'}}\\
& \le & c \n{J_-^{\frac{1}{2}}(J^su_1,J^su_3)}{L^2_{xt}}\n{J^{-3s}(J^su_2J^s\overline{u}_4)}{L_t^q(L_x^{2+})}
\end{eqnarray*}
with $\frac{1}{q}=-b'$, cf. the corresponding part of step 2. of this proof. Again we can use Corollary \ref{k23} and part ii) of Lemma \ref{l23} to obtain the desired bound.
$\hfill \Box$


\begin{thebibliography}{99}
\bibitem[BL]{BL}Bergh, J., L\"ofstr\"om, J.: Interpolation Spaces; Berlin - Heidelberg; 1976
\bibitem[BOP98]{BOP98}Bekiranov, D., Ogawa, T., Ponce, G.: Interaction equations for short and long dispersive waves, J. of Functional Analysis 158 (1998), 357 - 388
\bibitem[B93]{B93}Bourgain, J.: Fourier transform restriction phenomena for certain lattice subsets and applications to nonlinear evolution equations, GAFA 3 (1993), 107 - 156 and 209 - 262
\bibitem[B98]{B98}Bourgain, J.: Refinements of Strichartz' inequality and Applications to 2D-NLS with critical Nonlinearity, International Mathematics Research Notices 1998, No. 5, 253 - 283
\bibitem[CDKS01]{CDKS}Colliander, J. E., Delort, J. M., Kenig, C. E., Staffilani, G.: Bilinear estimates and applications to $2D$ NLS, Transactions of the AMS 353 (2001), 3307 - 3325
\bibitem[G96]{G96}Ginibre, J.: Le probl\`eme de Cauchy pour des EDP semi-lin\'eaires p\'eriodiques en variables d'espace (d'apr\`es Bourgain), Ast\'erisque 237 (1996), 163 - 187
\bibitem[GTV97]{GTV97}Ginibre, J., Tsutsumi, Y., Velo, G.: On the Cauchy-Problem for the Zakharov-System, J. of Functional Analysis 151 (1997), 384 - 436
\bibitem[Gr00]{G00}Gr\"unrock, A.: On the Cauchy- and periodic boundary value problem for a certain class of derivative nonlinear Schr\"odinger equations, preprint, arXiv:math.AP/0006195
\bibitem[KPV91]{KPV91}Kenig, C. E., Ponce, G., Vega, L.: Oscillatory Integrals and Regularity of Dispersive Equations, Indiana Univ. Math. J. 40 (1991), 33 - 69
\bibitem[KPV96]{KPV96}Kenig, C. E., Ponce, G., Vega, L.: Quadratic forms for the 1 - D semilinear Schr\"odinger equation, Transactions of the AMS 348 (1996), 3323 - 3353
\bibitem[St97]{St97} Staffilani, G.: Quadratic forms for a 2-D semilinear Schr\"odinger equation, Duke Math. J. 86 (1997), 79 - 107
\bibitem[T00]{T00}Tao, T.: Multilinear weighted convolution of $L^2$ functions, and applications to non-linear dispersive equations, preprint, arXiv:math.AP/0005001
\end{thebibliography}
\end{document}